\documentclass[12pt]{article}

\usepackage[latin1]{inputenc}
\usepackage{amssymb,amsmath, amsthm}

\def \({\left(}
\def \){\right)}
\def \[{\left[}
\def \]{\right]}
\def \R{\mathbb{R}}
\def \F{\mathcal{F}}
\def \O{\mathcal{O}}
\def \p{\partial}
\def \sign{\text{sign }}
\def \t{\tilde}
\def \ds{\displaystyle}

\theoremstyle{definition}
\newtheorem{thm}{Theorem}[section]
\newtheorem{prop}[thm]{Proposition}
\newtheorem{lemma}{Lemma}[section]
\newtheorem{rem}[lemma]{Remark}
\newtheorem{cor}[lemma]{Corollary}
\newtheorem{claim}[lemma]{Claim}
\newtheorem{defi}[lemma]{Definition}
\newtheorem*{pro}{Problem}

\begin{document}

\title{Weak uniqueness and partial regularity for the composite membrane problem}
\author{Sagun Chanillo$^\star$\\Department of Mathematics\\Rutgers University\\Piscataway, NJ
08854\\USA\\ \\Carlos E. Kenig\thanks{Supported in part by
NSF}\\Department of Mathematics\\University of Chicago\\Chicago, IL
60637\\USA}
\date{}
\maketitle

\noindent{\tt Feb6/2007.}
\section{Introduction}
Our main consideration will be the physical problem proposed in
\cite{CGIKO} which can be stated as:

\begin{pro}[P]\label{P} Build a body of prescribed shape out of given materials of
varying density, in such a way that the body has prescribed mass and
so that the basic frequency( with fixed boundary) is as small as
possible.
\end{pro}

This problem by virtue of Theorem 13 in \cite{CGIKO} can be
converted into the following minimization problem. Given a bounded
domain $\Omega \subset \R^n$ with smooth boundary, fix $\alpha>0$
and $A\in [0,|\Omega|]$. For any measurable subset $D\subset
\Omega$, denote by $\lambda_\Omega(\alpha,D)$ the first Dirichlet
eigenvalue for the problem,

\begin{equation}\label{eq0.1}  -\Delta u+\alpha
\chi_Du=\lambda_\Omega(\alpha, D)u, \ {\rm on}\ \Omega$$
$$u=0,\ {\rm on }\ \partial\Omega.\end{equation}
Define,

\begin{equation}\label{eq0.2} \Lambda_\Omega(\alpha,A)=\inf_{ D\subset
\Omega,\|D|=A}\lambda_\Omega(\alpha,D).\end{equation}

A minimizer $D$ to (\ref{eq0.2}) will be called an optimal
configuration for the data $(\Omega,\alpha,A)$. For this $D$ we
denote the associated eigenfunction solution to $(0.1)$ by $u$. The
pair $(u,D)$ will be called an optimal pair solution to the
composite problem or for short a solution to the composite problem.

A variational formulation of our problem is also possible and is
given by (see \cite{CGIKO}),

\begin{equation}\label{eq0.3} \Lambda_\Omega(\alpha,A)=\inf_{u\in H^1_0(\Omega), \|D|=A, \
||u||_2=1}\int_\Omega(|\nabla u|^2+\alpha\chi_D u^2).\end{equation}

Theorem 1 in \cite{CGIKO} establishes the basic properties of the
existence and regularity of optimal pairs.

\begin{thm}\label{thm0.1}
\cite{CGIKO} For any $\alpha>0$ and $A\in [0,|\Omega|]$, there
exists an optimal pair $(u,D)$. Moreover, the optimal pair $(u,D)$
has the property,
\begin{enumerate}
\item[(a)] $u\in C^{1,\gamma}(\overline \Omega)\cap
H^2(\overline \Omega)$, for every $\gamma<1$.
\item[(b)] $D$ is a sub-level set of $u$, that is there exists
$c\geq 0$ such that,
$$D=\{u\leq c\}.$$
\item[(c)] If $\alpha\not =\Lambda_\Omega(\alpha, A)$, then
every level set $\{u=s\}$ has measure zero.
\end{enumerate}
\end{thm}

See Remark \ref{rem1.2} for additional comments regarding (c). From
Theorem 13 in \cite{CGIKO} we also know that the physical problem
(P) stated earlier is equivalent to the variational problem
(\ref{eq0.3}) provided,
\begin{equation}\label{eq0.4} \alpha<\Lambda_\Omega(\alpha,
A).\end{equation} In the sequel we shall always assume
(\ref{eq0.4}). Now putting together Theorem \ref{thm0.1} and the
variational characterization of the problem $(0.3)$ we see that the
Euler-Lagrange equation of our problem is: $$-\Delta
u+\alpha\chi_{\{u\leq c\}}u=\Lambda_\Omega(\alpha, A)u,\ {\rm in}\
\Omega$$
\begin{equation}\label{eq0.5} u=0,\ {\rm on}\ \partial\Omega. \end{equation}

In Section 2 we first turn to the problem of uniqueness of optimal
pairs $(u,D)$. A principal result of \cite{CGIKO} is that even in
domains that exhibit symmetry, the optimal pair need not be unique,
and in fact uniqueness is known without any assumptions only if
$\Omega$ is the ball. Nevertheless, we establish that generically
there is a sort of weak uniqueness in the problem.

\begin{thm}[Weak Uniqueness]\label{thm0.2} Assume $(0.4)$. For almost every value of $A\in (0,|\Omega|)$, there exists $c>0$ such that for all
optimizing pairs $(u_i,D_i)$, $$D_i=\{ x|u_i(x)\leq c\}.$$
\end{thm}

Thus though there is non-uniqueness in the problem, the level height
where one must cut-off the eigenfunction to get $D_i$ must
generically be the same for all eigenfunctions.

Under additional assumptions, that is if eigenfunctions agree at one
point to infinite order or if $\Omega$ is convex in $\R^2$ with
additional assumptions, the assertion of weak uniqueness can be
turned into a statement of true uniqueness. See for example, Lemma
\ref{lemma1.19} and Theorem \ref{thm1.1} in Section 2.

In Section 3 we turn to the regularity of the free boundary ${\cal
F}$, defined by,

\begin{equation}\label{eq0.6} {\cal F}=\{ x| u(x)=c\}.\end{equation}

 We recall an initial result, Theorem $8$ in \cite{CGK},

\begin{thm}\label{thm0.3}
Let $x_0\in {\cal F}$. Assume $\nabla u(x_0)\not = 0$. That is $x_0$
is a regular point of the free boundary. Then there exists a ball
$B(x_0,r)$  of radius $r>0$ centered at $x_0$, and a real-analytic
function $\phi(x_1,x_2, \cdots , x_{n-1})$ such that,
$${\cal F}\cap B(x_0,r)=\{(x_1,x_2,\cdots , x_n)| x_n=\phi(x_1,x_2, \cdots ,x_{n-1})\}.$$
\end{thm}

 That is the free boundary in the neighborhood of a regular point is a
 hypersurface given by the graph of a real-analytic function.

Subsequently Blank \cite{B} performed a blow-up analysis in
dimension $2$ to classify the singular points of ${\cal F}$, that is
those points on ${\cal F}$ where $\nabla u=0$. This analysis in
dimension $2$ was completed in the paper by Shahgholian \cite{S},
who also obtained a condition that guarantees when singular points
of ${\cal F}$ in dimension $2$ are isolated.

The free boundary problem for the composite problem can be easily
converted to an equivalent problem (see for e.g. \cite{S}), given
by,
$$\Delta v=f\chi_{\{v\geq0\}}-g\chi_{\{v\leq0\}},$$
\begin{equation}\label{eq0.7}f,g\in C^{1,\gamma}, \ f>0,\ g<0,\
f+g<0.\end{equation} Our main result concerning the structure of
${\cal F}$ in Section 3 is:

\begin{thm}\label{thm0.4}[Structure of the Free Boundary of Solutions (\ref{eq0.7})]
For $\Omega \subset \R^n$, there is a decomposition,
$${\cal F}={\cal F}_0\cup S_v^1\cup S_v^2,$$
where $S_v^2$ has Hausdorff $\dim \leq n-2$, ${\cal
H}^{n-1}(S_v^1)\leq C$, and for all $x_0\in {\cal F}_0$, there
exists a ball $B(x_0,r)$ centered at $x_0$ such that, ${\cal F}\cap
B(x_0,r)$ is a hypersurface given by the graph of a real-analytic
function.
\end{thm}

The principal tool we use to perform our blow-up analysis and
thereby get Theorem \ref{thm0.4} is an energy functional introduced
by Weiss \cite{W}. Set, ($f\equiv f_0,\,g\equiv g_0$)

\begin{equation}\label{eq0.8}
W(r)={{1}\over{r^{n+2}}}\int_{B(x_0,r)}\(|\nabla
v|^2+2(f_0v^++g_0v^-)\) -{{2}\over{r^{n+3}}}\int_{\partial
B(x_0,r)}u^2.\end{equation}

Weiss showed that $W(r)$ is monotonically increasing. We offer an
alternative proof based in part on the Rellich-Pohozhaev identity
which explicitly shows that no structural assumptions are needed to
get the monotonicity.

Next we proceed to classify the blow-up limits in the spirit of the
paper by Monneau-Weiss \cite{MW}. Two points are to be noted in
contrast to \cite{MW}. First that in our case blow-up limits are
non-degenerate, and second that we have two types of blow-up limit
solutions that are homogeneous of degree $2$. This is already
evident in the work in dimension $2$ by Blank \cite{B} and
Shahgholian \cite{S}.

Lastly we address the question of $C^{1,1}$ bounds. In general such
bounds are not available for the composite problem if we only
analyze the Euler-Lagrange equation (\ref{eq0.7}). So-called cross
solutions arise from homogeneous harmonic polynomials of degree $2$
with corresponding failure of $C^{1,1}$ bounds in dimension $2$ as
has been exhibited by Andersson and Weiss \cite{AW} in the case
$f\equiv-1,\,g\equiv0$. The example of Andersson-Weiss can be easily
extended to all dimensions by the addition of dummy variables. We
show that the \cite{AW} construction extends to our setting (Remark
\ref{rem2.40}). Our regularity result proved in Section 3 is:

\begin{thm}\label{thm0.5}
The free boundary $\F=G\cup B,$ where in $G$ we have pointwise
$C^{1,1}$ bounds and $B$ has Hausdorff $\dim\leq n-2.$ (See Theorem
\ref{thm2.38}, Definition \ref{defi2.35} and Remark \ref{rem2.39}).
\end{thm}

It remains open whether proceeding from the variational problem
instead of (\ref{eq0.7}) allows one to get $C^{1,1}$ bounds. It is
readily seen that global assumptions on the boundary of $\Omega$ do
ensure that $C^{1,1}$ bounds and full regularity are achieved. A
result of this type proved in Section 3 is, (Proposition
\ref{prop2.46})

\begin{prop}\label{prop0.6} Assume $\Omega\subset\R^2$ has
two axes of symmetry. Then the free boundary ${\cal F}$ is a
real-analytic curve and $u\in C^{1,1}$.
\end{prop}

\setcounter{equation}{0}

\section{Uniqueness and Weak uniqueness}
Our goal is to prove Theorem \ref{thm0.2} of the introduction in
this section. We shall also show that a weak uniqueness assertion
like in Theorem \ref{thm0.2} can be converted to a uniqueness
assertion on convex domains with additional assumptions. Let
$\Omega\subset\R^n$, be a bounded domain with $\partial\Omega$
smooth. For $\alpha>0,\ A\in \[0,|\Omega|\]$, and $D\subset\Omega$,
let $\lambda_{\Omega}(\alpha, D)=\lambda$, be the lowest eigenvalue
to,

\begin{equation}
\begin{cases}\label{eq1.1}
-\Delta v+\alpha\chi_Dv=\lambda v & \text{ on } \Omega\\
v|_{\partial\Omega}=0
\end{cases}
\end{equation}

The variational characterization of (\ref{eq1.1}) gives,

\begin{equation}\label{eq1.2}
\lambda=\inf_{u\in
H_0^1(\Omega)}\frac{\displaystyle\int_{\Omega}\(|\nabla
u|^2+\alpha\chi_Du^2\)}{\displaystyle\int_{\Omega}u^2}
\end{equation}

\begin{lemma}\label{lemma1.1} There exists a unique minimizer
$v\in H_0^1$ of (\ref{eq1.2}) with $\|v\|_2=1$, which is
non-negative.
\end{lemma}

\begin{proof}
By Theorem 8.38 in \cite{GT}, the eigenvalue $\lambda$ is simple and
the eigenspace is spanned by a non-negative eigenfunction. Since
$\|v\|_2=1$, we have a unique non-negative eigenfunction, with the
property, $\|v\|_2=1$.
\end{proof}

Define
$$\Lambda=\Lambda_{\Omega}(\alpha, A)=\inf_{\substack{D\subset\Omega\\|D|=A}}\lambda(\alpha,A).$$

\begin{rem}\label{rem1.2}
Assume $\alpha<\Lambda.$ Then for the solution to the composite
problem $(u, D)$ stated in the introduction, $|\{u=s\}|=0$, for all
$s$. This is Theorem  1(c) in \cite{CGIKO}. Note that $s=0$ is not
covered by the proof in \cite{CGIKO} but is easily ruled out by
superharmonicity of $u$.
\end{rem}

\begin{lemma}\label{lemma1.3}
Let $\F=\{u:u=c\}$, where $D=\{x\in\Omega:u\leq c\}$ and $(u,D)$ is
the solution of our composite problem. Then $\nabla u\not\equiv0$ on
$\F$. (In fact in the boundary of each connected component of $^CD$,
$\nabla u$ cannot be identically 0).
\end{lemma}

\begin{proof}
Assume $\nabla u|_{\F}\equiv 0.$ Consider the open set
$\O=\{x:u>c\}$. Since $A<|\Omega|$, by remark \ref{rem1.2},
$|\O|>0$. Let $U$ be a connected component of $\O$, and let
\begin{equation}\label{eq1.3}
\begin{cases}
-\Delta w=\mu_U w& \text{ in } U\\
w|_{\partial U}=0.
\end{cases}
\end{equation}
where $\mu_U$ is the first Dirichlet eigenvalue of $U$. We claim
$\Lambda\leq\mu_U.$ To check this extend $w$ to $U^c$, by setting
$w\equiv 0$ in $U^c$. The extended function will still be denoted by
$w$ and we may normalize it so that $\|w\|_2=1.$ Then,
$$\Lambda\leq\int_{\Omega}|\nabla w|^2+\alpha\int_Dw^2=\int_{\Omega}|\nabla w|^2=\mu_U.$$
If $\Lambda=\mu_U$, by Lemma \ref{lemma1.1}, $u=w$. Since $w\equiv0$
on $D$, and since $u$ is superharmonic, so is $w$, so $w=u=0$, a
contradiction. So, $\Lambda<\mu_U$. Let $v=\partial_{x_j}u$, for
some fixed $j.$

In $U$,

\begin{equation}\label{eq1.4}
-\Delta u=\Lambda u,
\end{equation}

so that on differentiating (\ref{eq1.4}), $v$ satisfies,

\begin{equation}\label{eq1.5}
\begin{cases}
-\Delta v=\Lambda v & \text{ in } U\\
v|_{\partial U}\equiv0 & v\in C^{\gamma}(\bar U)
\end{cases}
\end{equation}

We claim $v\equiv0.$ This will imply $u\equiv c$ in $U$, which will
contradict Remark \ref{rem1.2}. Since $\Lambda<\mu_U$, we may solve
using the Fredholm alternative,

\begin{equation}\label{eq1.6}
-\Delta f-\Lambda f=-\Lambda\text{ in $U$, } f\in H_0^1(U)
\end{equation}

Let $h=f^+=\max(f,0)$. Clearly $h\in H_0^1(U).$ Multiplying
(\ref{eq1.6}) by $h$ and integrating by parts,
$$\int_U|\nabla h|^2-\Lambda\int_U h^2=-\int_U\Lambda h\leq0.$$
Thus,

\begin{equation}\label{eq1.7}
\int_U|\nabla h|^2\leq\Lambda\int_U h^2.
\end{equation}

If $\displaystyle\int_Uh^2\neq0$, then from (\ref{eq1.7}),
$\mu_U\leq\Lambda$. This is a contradiction. Hence
$\displaystyle\int_Uh^2=0$ and $h\equiv0$ in $U.$ Thus $f\leq0.$ Set
$\psi=1-f$. Thus $\psi\geq1$, and from (\ref{eq1.6}),
$$-\Delta\psi-\Lambda\psi=0.$$
By elliptic regularity, $\psi\in C^{\infty}(U).$ Now find
$U_j\Subset U$, with $dist(\partial U_j,
\partial U)\to0$, and $\partial U_j$ smooth. So if $x\in U$, then $x\in
U_j$ for large enough $j.$ Let $\phi=v/\psi$, where $v$ is defined
in (\ref{eq1.5}). Note, because $\psi\geq1$, and by (\ref{eq1.5})
again,

\begin{equation}\label{eq1.8}
\sup_{\partial U_j}|\phi|\leq\sup_{\partial U_j}|v|\to0 \text{ as
}j\to \infty
\end{equation}

Now,
$$\nabla\phi=\frac{\nabla v}{\psi}-\frac{v\nabla\psi}{\psi^2}=\frac{\psi\nabla v-v\nabla\psi}{\psi^2}$$
and
$$\Delta\phi=\frac{\nabla\psi\cdot\nabla v+\psi\Delta v-\nabla v\cdot\nabla\psi-v\Delta\psi}{\psi^2}
-\frac2{\psi^3}\(\psi\nabla v-v\nabla\psi\)\cdot\nabla\psi$$
$$=\frac{\psi\Delta v-v\Delta\psi}{\psi^2}-\frac2{\psi}\nabla\psi\cdot\nabla\phi=
-\frac2{\psi}\nabla\psi\cdot\nabla\phi.$$

Thus $\phi$ satisfies
$$\Delta\phi+\frac2{\psi}\nabla\psi\cdot\nabla\phi=0\text{ in }U_j$$

Thus by the maximum principle, and (\ref{eq1.8}),
$$\sup_{\bar U_j}|\phi|\leq\sup_{\partial U_j}|\phi|\to0\text{ as }j\to\infty$$

Thus $\phi\equiv0$ in $U$ and hence $v\equiv0$ in $U$. Our
proposition is proved.
\end{proof}

Combining Theorem 8 in \cite{CGK} and Lemma \ref{lemma1.3} , we
have,

\begin{lemma}\label{lemma1.4}
If $(u,D)$ is a minimizing pair with $\alpha<\Lambda_{\Omega}$, then
there exists $x_0\in\F=\{u=c\}$ and a ball $B(x_0, r)=B$ centered at
$x_0$, so that $B\subset\Omega$ and
$$\F\cap B=\{(x,\phi(x)), x\in\R^{n-1}, \phi:U\subset\R^{n-1}\to\R\},$$
with $\phi(x_0)=0, \nabla\phi(x_0)=0$ and $\phi(x)$ real-analytic.
Furthermore,
$$D\cap B=\{(x,y):y<\phi(x)\}\cap B$$
$$^cD\cap B=\{(x,y):y>\phi(x)\}\cap B$$
\end{lemma}

\begin{lemma}\label{lemma1.5}
Let $\psi(x'):U\subset\R^{n-1}\to\R$ be smooth; where $U$ is open
and $U\supset B(0,r).$ Assume $\psi(0)=0,\, \nabla\psi(0)=0$, and
let
$$D=\{(x',y):y<\psi(x')\}\cap B,\ B=B(0,r).$$
Then there exists $\epsilon_0>0$, and a smooth function, $x=(x',y)$
$$\Phi(t,x)=\Phi_t(x):\{|t|\leq\epsilon_0\}\times B\to B$$
such that,
\begin{enumerate}
\item[(a)] for all fixed $t$,
$$\Phi_t:\overline B\to\overline B$$
is a diffeomorphism, with $\Phi(0,x)=\Phi_0(x)=x.$
\item[(b)] for all $t$, $|t|\leq\epsilon_0$, and some $\delta>0$, $\delta<r/50$
$$\Phi_t|_{\overline{B}\setminus B(0,2\delta)}=x.$$
\item[(c)] Let $\chi_(D_t)(x)=\chi_D\(\Phi_{-t}(x)\).$ Then
$$\frac{d}{dt}\(|D_t|\)|_{t=0}=1.$$
\end{enumerate}
\end{lemma}

\begin{proof}
Let $f(x)$ be a smooth cut-off function, $f\in
C_0^{\infty}\(B(0,\delta/100)\),\,f\geq0$. Let $\nu(x')$ denote the
unit outward normal to $y=\psi(x').$ We extend $\nu(x')$ smoothly as
a vector field $X$ to all points in $B(0,\delta/10)$. Now define
\begin{equation}\label{eq1.9}\displaystyle\frac{d\Phi_t}{dt}(x)=\frac{X(x)f(x)}{\int\limits_{\partial
D\cap B(0,\delta/10)}f(\sigma)\ d\sigma}=V(x),\
\Phi_0(x)=x\end{equation}

(a), (b) follows from (\ref{eq1.9}). Note that a simple degree
argument is needed to show that $\Phi_t$ is a diffeomorphism. (c)
follows from the Appendix 1, by noting that
$$V|_{\partial D}=\frac{\nu(x')f(x)}{\int\limits_{\partial
D}f(\sigma)\ d\sigma}.$$ Hence
$$\int_{\partial D}\langle V,\,\nu\rangle=1.$$
\end{proof}

\begin{lemma}\label{lemma1.6}
Construct $\Phi_t(x)$ as in Lemma \ref{lemma1.5}, $x_0=0$ as in
Lemma \ref{lemma1.4}. Define
$$\phi_t(x):\Omega\to\Omega,\text{ by }$$
$$\phi_t(x)=\begin{cases}
\Phi_t(x)& x\in B(0,3\delta)\\
x& x\in \Omega\setminus B(0,3\delta).
\end{cases}$$
\begin{enumerate}
\item[(a)] Then $\phi_t(x)$ is a diffeomorphism of $\Omega$.
\item[(b)] If $D_t=\{\phi_t(x),\, x\in D\}$, then
$$\frac{d}{dt}\(|D_t|\)=1.$$
\item[(c)] If $(u,D)$ is a solution to the composite problem and
$$-\Delta u_t+\alpha\chi_{D_t}u_t=\lambda(t)u_t$$
$$u_t|_{\partial\Omega}=0,$$
where $D=\{x\in\Omega,\,u\leq c\}=D_0,\,u_0=u,\,\lambda(0)=\Lambda.$
Then
$$\lambda'(0)=\alpha c^2.$$
\end{enumerate}
\end{lemma}

\begin{proof}
Using (b) in (A1.10) and Lemma \ref{lemma1.5}(c) we get (c). (b)
follows from Lemma \ref{lemma1.5}. (a) follows from the definition
of $\phi_t(x)$ and Lemma \ref{lemma1.5}(a).
\end{proof}

\begin{lemma}\label{lemma1.7}
Assume that $\Lambda_{\Omega}(\alpha,A)$ is differentiable at
$A=A_0.$ Let $(u,D)$ be a minimizer. Construct domains $D_t$ as in
Lemma \ref{lemma1.6}, where $B=B(x_o,r)$ is supplied by Lemma
\ref{lemma1.4}. Then
$$\frac{d}{dt}\Lambda(\alpha,A)|_{A=A_0}=\alpha c^2.$$
\end{lemma}

\begin{proof}
Let $|D_t|=m(t)$; $D_t$ as in Lemma \ref{lemma1.6}. Let
$f(t)=\Lambda\(\alpha,m(t)\).$ Then, $f$ is differentiable at $t=0.$
\begin{equation}\label{eq1.10}
f'(0)=\left.\frac{d\Lambda}{dA}(\alpha,A)\right|_{A=A_0}\cdot
m'(0)=\left.\frac{d\Lambda}{dA}(\alpha,A)\right|_{A=A_0}.\end{equation}
Next for $t>0$, by the definition of $\Lambda,$
$$\frac{f(t)-f(0)}t\leq\frac{\lambda(t)-f(0)}t=\frac{\lambda(t)-\lambda(0)}t.$$
Letting $t\downarrow0$, we get $f'(0)\leq\lambda'(0).$ Arguing
similarly for $t<0$, letting $t\uparrow0$, using the
differentiability at $t=0$ of $f$ and $\lambda$ we get
$f'(0)=\lambda'(0)=\alpha c^2$ by Lemma \ref{lemma1.6}. Thus from
(\ref{eq1.10}),
$$\left.\frac{d\Lambda}{dA}(\alpha,A)\right|_{A=A_0}=\alpha c^2.$$
\end{proof}

\begin{proof}[Proof of Theorem \ref{thm0.2}]
$\Lambda(\alpha, A)$ is strictly increasing in $A$ and Lipschitz in
$A$, Prop. 10, \cite{CGIKO}. Thus $\Lambda'(\alpha,A)$ exists, a.e.
A., and $\Lambda'(\alpha, A)=\alpha c^2$ by Lemma \ref{lemma1.7}.
Hence if $(u_1, D_1),\,(u_2, D_2)$ are two configurations $|D_i|=A$,
with $D_i=\{x:u_i<c_i\},$ then $\alpha c_1^2=\alpha c_2^2.$ Hence
$c_1=c_2.$
\end{proof}

We shall now show that under some conditions, the weak uniqueness
conclusion of Theorem \ref{thm0.2} can be turned into a uniqueness
result. We will restrict our attention to domains
$\Omega\subset\R^2.$

\begin{lemma}\label{lemma1.8}
Let $\Omega\subset\R^2$, and let
\begin{equation}\label{eq1.11}\begin{cases}
-\Delta u+\alpha\chi_{\{u\leq c\}}(x)u=\lambda u\\
u|_{\partial\Omega}=0,\ \|u\|_2=1.
\end{cases}\end{equation}
Then for any $x_0\in\R^2$,
$$\frac12\int_{\partial\Omega}\langle x-x_0,\,\nu\rangle\(\frac{\p u}{\p \nu}\)^2=\lambda-\alpha c^2|D^c|-\alpha\int_D u^2$$
where $D=\{x:u\leq c\}.$
\end{lemma}

\begin{proof}
We use the Rellich-Pohozhaev identity. Now,
$$-\langle x-x_0,\nabla u\rangle\Delta u=-\nabla\cdot\(\langle x-x_0,\nabla u\rangle\nabla u\)+|\nabla u|^2+\frac12\(x-x_0\)\cdot\nabla\(|\nabla u|^2\).$$
Thus, integrating the identity above over $\Omega$,
\begin{multline}\label{eq1.12}
-\int_{\Omega}\langle x-x_0,\nabla u\rangle\Delta u\\
=-\int_{\Omega}\langle x-x_0,\nu\rangle\(\frac{\p u}{\p \nu}\)^2\
d\sigma+\frac12\int_{\Omega}\langle x-x_0,\nu\rangle\(\frac{\p
u}{\p\nu}\)^2\ d\sigma\\=\frac12\int_{\Omega}\langle
x-x_0,\nu\rangle\(\frac{\p u}{\p\nu}\)^2\ d\sigma
\end{multline}
From (\ref{eq1.11}),
\begin{equation}\label{eq1.13}
-\Delta u=\lambda u-\alpha\chi_D u.
\end{equation}
Substituting (\ref{eq1.13}) into the left side of (\ref{eq1.12}) we
get,
$$\int_{\Omega}\alpha\chi_D u\langle x-x_0,\nabla u\rangle-\int_{\Omega}\langle x-x_0,\nabla u\rangle\lambda u
=\int_{\Omega}\langle x-x_0,\nu\rangle\(\frac{\p u}{\p\nu}\)^2.$$
Thus,
\begin{multline}\label{eq1.14}
\frac12\int_{\Omega}\langle x-x_0,\nu\rangle\(\frac{\p
u}{\p\nu}\)^2\\=-\frac{\lambda}2\int_{\Omega}\langle
x-x_0,\nabla\(u^2\)\rangle+\frac{\alpha}2\int_{\Omega}\langle
x-x_0,\nabla\(u^2\)\rangle.
\end{multline}
The first integral on the right by integration by parts is
\begin{equation}\label{eq1.15}
\lambda\int_{\Omega}u^2=\lambda.
\end{equation}
For the second integral, since $D=\{x:u(x)\leq c\},$ there exists
$c_j\uparrow c$, such that by Sard's Theorem $c_j$ is a regular
value. Let $D_j=\{x:u<c_j\}$. Now by integration by parts,
\begin{eqnarray*}
\int_{\Omega}\langle
x-x_0,\nabla\(u^2\)\rangle&=&-2\int_{D_j}u^2+\int_{\p
D_j\cap\Omega}\langle x-x_0,\nu\rangle u^2\\
&=&-2\int_{D_j}u^2+\int_{\p
D_j\cap\Omega}c_j^2\langle x-x_0,\nu\rangle\\
&=&-2\int_{D_j}u^2+c_j^2|^cD_j|.
\end{eqnarray*}
Letting $j\to\infty,$
\begin{equation}\label{eq1.16}\int_D<x-x_0,\nabla(u^2)>=\ -2\int_Du^2+c^2|^cD|\end{equation}
Inserting (\ref{eq1.16}),\ (\ref{eq1.15}) into (\ref{eq1.14}) we get
our result.
\end{proof}

To obtain a true uniqueness assertion we first need a preliminary
lemma which is valid in all dimensions. We shall assume in the
sequel that our solutions are normalized by the condition
$\|u\|_2=1$.

\begin{lemma}\label{lemma1.19}
Let $(u_i,D_i),\ i=1,2$ be two solutions of our composite problem.
Assume that $D_1$ is connected. Assume furthermore we have weak
uniqueness that is $D_i=\{x\in \Omega: u_i\leq c\}$ and $u_1-u_2$
vanishes at a single point $x_0\in D_1$ to infinite order. Then
$u_1\equiv u_2$ in $\Omega$.
\end{lemma}

\begin{proof}
First we note $u_1(x_0)=u_2(x_0)<c$. Thus there is a ball $B$
centered at $x_0$ where $u_i(x)<c,\ i=1,2$. Thus in this ball we
have,
\begin{equation}\label{eq1.17}
-\Delta u_i+\alpha u_i=\Lambda u_i,\ i=1,2\end{equation}
 Thus,
$w=u_1-u_2$ also satisfies the equation (\ref{eq1.17}) and $w$
vanishes at $x_0$ to infinite order. Thus, $w$ vanishes identically
in $B$. Now consider the set,
$$W={\rm int}\{ x\in D_1, \ u_1=u_2\}.$$
We have established that $W$ is non-empty. We shall now show that
$W$ is both open and closed in the relative topology of $D_1$. Since
$D_1$ is connected we then get $W=D_1$. Since  $u_1=u_2<c$ on
$\stackrel{\circ}{D_1}$ we obtain that, $D_1\subset D_2$. Since
$|D_1|=|D_2|$ we see right away that $D_1=D_2$.

Now by definition $W$ is open. So let $z_0\in \overline W=F\cap D_1$
where $F$ is closed. Thus $u_1(z_0)=u_2(z_0)<c$ and thus there is a
ball $B$ centered at $z_0$ where (\ref{eq1.17}) is satisfied. Again
$w$ satisfies (\ref{eq1.17}) with $w$ vanishing on some open set in
$B$. This is because $z_0$ is a boundary point to $W$. Again by
unique continuation $w$ vanishes in $B$. Thus $z_0\in W$. We have
checked $W$ is also closed. Since now $D_1=D_2$, applying Lemma
\ref{lemma1.1} we obtain the conclusion of our lemma.
\end{proof}

\begin{rem}\label{rem1.10}
The same result holds if $x_0\in\partial\Omega.$ The proof is
similar, but slightly more complicated.
\end{rem}

\begin{thm}\label{thm1.1}
Assume $\Omega\subset\R^2$ with smooth boundary. Assume that
$\Omega$ is strictly convex. Let $(u_i,D_i)$ be two solutions to the
composite problem with eigenvalue $\Lambda$. Assume that,
\begin{enumerate}
\item[(a)] $$\int_{D_1}u_1^2=\int_{D_2}u_2^2.$$
\item[(b)] Weak uniqueness holds, $D_i=\{x\in \Omega| u_i(x)\leq c\}$.
\item[(c)] The sets $\{x| u_1(x)<u_2(x)\}$ and $\{x|u_1(x)>u_2(x)\}$
are both connected.
\end{enumerate}
Then $u_1\equiv u_2$.
\end{thm}

\begin{proof}
Since $\Omega$ is convex, it is simply connected and since
$\alpha<\Lambda$, by Theorem 2 \cite{CGIKO} the sets $D_i$ are
connected. Writing Lemma \ref{lemma1.8} for $u_i$ and subtracting
the expression for $u_2$ from that of $u_1$, we get after using the
hypothesis $(a), \ (b)$ above that,
$$\int_{\partial \Omega} \langle x-x_0,\nu\rangle \[ \({{\p
u_1}\over{\partial \nu}}\)^2-\({{\p u_2}\over{\p \nu}}\)^2\]=0.$$ We
re-write this expression to get,
\begin{equation}\label{eq1.18}\int_{\partial \Omega} \langle x-x_0,\nu\rangle{{\p}\over{\p
\nu}}(u_1+u_2) {{\p}\over{\p\nu}}(u_1-u_2)=0.\end{equation}

Now in a tubular neighborhood of $\partial\Omega$ both $u_1, u_2$
satisfy (\ref{eq1.17}). Thus $u_1+u_2$ also satisfies (\ref{eq1.17})
with $u_1+u_2>0$ in $\Omega$ and vanishing on $\partial\Omega$. Thus
by Hopf's boundary point lemma,
\begin{equation}\label{eq1.19}{{\p}\over{\p\nu}}(u_1+u_2)<0\end{equation}

Now set $\psi=u_1-u_2$. Let, $$E_1=\left\{x\in \p\Omega\left| {{\p
\psi}\over{\p\nu}}>0\right.\right\},\ E_2=\left\{ x\in\p\Omega\left|
{{\p \psi}\over{\p\nu}}<0\right.\right\}.$$ We show both sets are
empty. If we establish this result we have the conclusion of the
lemma. The reason is that if ${{\partial\psi}\over{\p v}}=0$ on
$\partial\Omega$, since $\psi=0$ on $\partial\Omega$ we conclude by
the Cauchy-Kovalevskaya Theorem that $\psi$ vanishes in a
neighborhood of a boundary point and thus applying Lemma
\ref{lemma1.19} we conclude $u_1=u_2$ in $\Omega$.

\paragraph{Case 1:} Assume w.l.o.g. that $E_2$ is empty and $E_1$ is
non-empty. Pick any $x_0\in \Omega$. Then by the strict convexity of
$\partial\Omega$, $\langle x-x_0,\nu\rangle>0$. Thus by
(\ref{eq1.19}) and the choice of $x_0$ we conclude that the integral
in (\ref{eq1.18}) is negative. This contradicts the identity
(\ref{eq1.18}).

\paragraph{Case 2:} We may now assume that both $E_1$ and $E_2$ are
non-empty. Consider the components of $E_1$ and $E_2$ on
$\partial\Omega$. These are intervals. We claim that the hypothesis
(c) rules out interlacing of intervals. That is the intervals that
make up the components of $E_1$ must share at least one boundary
point and likewise for the intervals that make up the components of
$E_2$. For assume there exist two intervals $I_1, I_2$ which are
components of $E_1$ and two intervals $J_1, J_2$ which are
components of $E_2$. Now we shall obtain a contradiction if we
assume that $I_1, I_2$ lie in different components of
$\partial\Omega\setminus \(J_1\cup J_2\)$. Taking  interior points
in $I_1, I_2$ we can connect the points by a curve that lies
entirely in $\Omega$ and in the set $\{u_1<u_2\}$. Now it is easily
seen that $\{u_1>u_2\}$ is disconnected. This contradicts (c). Thus
we have shown that $\partial\Omega$ consists of two arcs $\gamma_1,\
\gamma_2$ such that, $\gamma_1$ and $\gamma_2$ have common endpoints
$P,Q$ and such that, on $\gamma_1$, ${{\partial\psi}\over{\p
\nu}}\geq 0$, with ${{\partial\psi}\over{\p \nu}}>0$ on some
sub-interval of $\gamma_1$. Likewise,on $\gamma_2$,
${{\partial\psi}\over{\p \nu}}\leq 0$, with ${{\partial\psi}\over{\p
\nu}}<0$ on some sub-interval of $\gamma_2$. Now consider tangent
lines to $\partial\Omega$ at $P,Q$.

If the tangent lines intersect at $x_0$, apply (\ref{eq1.18}) with
this choice $x_0$. Notice then by the strict convexity of
$\partial\Omega$, $\langle x-x_0,\nu\rangle>0$ (except possibly at
$P,Q$) on $\gamma_1$ and $\langle x-x_0,\nu\rangle<0$ on $\gamma_2$.
Thus using (\ref{eq1.19}) and the behavior of $\psi$ on $\gamma_1,\
\gamma_2$ we easily see that the integral in (\ref{eq1.18}) is
negative. This is a contradiction.

Assume thus that the tangent lines at $P,Q$ are parallel and with no
loss of generality assume they are parallel to the $x_1$ axis,
$x=(x_1,x_2)$. Set $v(x)=\(n_1(x),n_2(x)\)$. Now (\ref{eq1.18})
holds for every $x_0$. Set $x_0=\(x^0_1, x^0_2\)$.  We may now
differentiate (\ref{eq1.18}) with respect to $x^0_1$ and we obtain,
$$\int_{\partial \Omega}n_1(x){{\partial}\over{\partial \nu}}(u_1+u_2){{\partial\psi}\over{\partial
v}}=0.$$ Now we may assume that $n_1(x)>0$ on $\gamma_1$ and
$n_1(x)<0$ on $\gamma_2$ except at $P,Q$ by the strict convexity of
$\partial\Omega$. Thus the integrand in the integral in
(\ref{eq1.18}) is non-positive by the use of (\ref{eq1.19}).
Furthermore, from (\ref{eq1.19}) and the behavior of $\psi$ on the
arcs $\gamma_i$ there are arcs on $\partial\Omega$ where the
integrand in (\ref{eq1.18}) is negative. This again is a
contradiction to (\ref{eq1.18}). Thus both sets $E_1$ and $E_2$ are
empty. Our Theorem is established.
\end{proof}

\setcounter{equation}{0}
\section{Partial Regularity}

Our goal is to prove Theorems \ref{thm0.4} and \ref{thm0.5} of the
introduction in this section. We follow the works of Blank \cite{B},
Shahgolian \cite{S}, Weiss \cite{W} and Monneau-Weiss \cite{MW},
with some necessary variants and extensions.

\paragraph{The set-up:} Let $\Omega\subset\R^n,$ be a bounded domain
with $\partial\Omega$ smooth. For $\alpha>0,\, A\in\(0,|\Omega|\)$,
we let $(u,D)$ be a solution of the composite problem, so that
\begin{equation}\label{eq2.1}
\begin{cases}
-\Delta u+\alpha\chi_{\{u\leq c\}}u=\Lambda u \text{ on }\Omega\\
u|_{\partial\Omega}=0,\ \int_{\Omega}u^2=1
\end{cases}
\end{equation}
where $D=\{u\leq c\}$. Recall that $u\geq0$ in $\Omega$ and that we
are assuming throughout that $\alpha<\Lambda.$ Note that $u\in
W^{2,p}(\Omega)\,\forall\,1\leq p<\infty,\, u\in
C^{1,\gamma}\(\overline\Omega\),\, 0\leq\gamma<1,$ with norm
depending only on $A,\,n,\,\Omega,\,p,\,\gamma,\,\alpha$ and
$\Lambda.$ Note also that $c>0$ since if $u(x_0)=0$, by
superharmonicity of $u$, $x_0\in\partial\Omega$, and $|\{u\leq
c\}|=A>0.$ Note also that $|\{u=c\}|=0$, by Remark \ref{rem1.2}. We
next let $v=c-u$ and write the equation for $v$, namely
\begin{equation}\label{eq2.2}
\Delta v=f\chi_{\{v\geq0\}}-g\chi_{\{v<0\}}
\end{equation}
where $f=(\Lambda-\alpha)u,\,g=-\Lambda u.$ Fix a neighborhood $U$
of $\F=\{u=c\}$, the free boundary, so that $f>0,\,g<0$ and $f+g<0$
in $\overline U.$ We thus have a solution $v$ of (\ref{eq2.2}), in
$U$ open, and functions $f,\,g\in C^{1,\gamma}\(\overline U\),$ with
norm bounded by $\tilde B_1=\tilde B_1(\gamma, u, \alpha, \Lambda,
A, \Omega)$ in $\overline U$, with $f,\,g\in W^{2,p}(U)$, with norm
bounded by $\tilde B_2=\tilde B_2(p, n, \alpha, \Lambda, A, \Omega)$
and with $|\Delta f|,|\Delta g|$ bounded by $\tilde B_3=\tilde
B_3(\alpha, \Lambda)$, and such that, for some
$\eta_0>0,\,\eta_0=\eta_0(\alpha, \Lambda, A, n, \Omega),$ we have
$f\geq\eta_0>0,\,g\leq-\eta_0,\,(f+g)\leq\eta_0$ in $\overline U.$
We also have $\|v\|_{C^{1,\gamma}(\overline
U)}+\|v\|_{W^{2,p}(\overline U)}\leq
N,\,N=N(\gamma,\,p,\,n,\,x,\,\Lambda,\,A,\,\Omega).$ Finally, we fix
$r_0$ so small that for all $x_0\in\F$ we have that $B\(x_0,
r_0\)\subset U$. We still study the behavior of $S_u=\{x\in\F:\nabla
u(x)=0\}=S_v=\{x\in\F:\nabla v(x)\}=0$, where $\F=\{v=0\}.$ Note
that, by \cite{CGK} (see Lemma 8 and Theorem (0.3) here) for each
$x_0\in\F\setminus S_v$, there exists a neighborhood $V_{x_0}$
around $x_0$ so that $\F$ is real analytic in it and $v$ (and $u$)
are real analytic in $V_{x_0}\cap\overline
D,\,V_{x_0}\cap\overline{^cD}$. One of our main tools in this
section is an energy functional introduced by Weiss:
\begin{equation}\label{eq2.3}
W(r)=\frac1{r^{n+2}}\int_{B\(x_0, r\)}\(|\nabla
v|^2+2\(fv+gv^-\)\)-\frac2{r^{n+3}}\int_{\partial B\(x_0, r\)}v^2
\end{equation}

In the next Lemma we compute $W'(r)$. (See \cite{W}, where the
computation is also carried out).

\begin{lemma}\label{lemma2.4}
Let $x_0\in S_{v},\, 0<r<r_0.$ Then, for $0<r<r_0,$
\begin{equation}\label{eq2.5}
W'(r)=\frac2{r^{n+2}}\int_{\partial B_r}\[\frac{\p
v}{\p\nu}-2\frac{v}r\]^2\ d\sigma+e(r),
\end{equation}
where for $0\leq\gamma<1$ we have, for $0<r<r_0$,
\begin{equation}\label{eq2.6}
|e(r)|\leq F\(n,\,\gamma,\,\|\nabla f\|_{\infty},\,\|\nabla
g\|_{\infty},\,N\)r^{\gamma-1},
\end{equation}
with $F(-,\,-,\,0,\,0,\,-)\equiv 0.$ (Here $\nu$ is the outward unit
normal to $\partial B_r$ and $B_r$ stands for $B(x_0,\,r)$).
\end{lemma}

\begin{proof}
We can assume, without loss of generality, that $x_0=0.$ We have:
$$\frac{\partial}{\partial r}\(\frac1{r^{n+2}}\int_{B_r}|\nabla v|^2\)=-\frac{n-2}{r^{n+3}}\int_{B_r}|\nabla v|^2+\frac1{r^{n+2}}\int_{\partial B_r}|\nabla v|^2.$$
Moreover, the Rellich-Pohozaev identity gives:
$$div\(x|\nabla v|^2\)=2div\(x\cdot\nabla v\nabla v\)+(n-2)|\nabla v|^2-2x\cdot\nabla v\Delta v,$$
we also have the identity
$$\(f\chi_{\{v\geq0\}}-g\chi_{\{v<0\}}\)\nabla v=\nabla\(fv^++gv^-\)-\nabla fv^+-\nabla gv^-,$$
$$\int_{B_r}x\cdot\nabla\(fv^++gv^-\)=r\int_{\partial B_r}\(fv^++gv^-\)-n\int_{B_r}\(fv^++gv^-\),$$
so that
\begin{multline*}\int_{\partial
B_r}|\nabla v|^2=2\int_{\partial B_r}\(\frac{\p
v}{\p\nu}\)^2+\frac{n-2}r\int_{\partial
B_r}|\nabla v|^2-2\int_{B_r}\(fv^++gv^-\)\\
+\frac{2n}r\int_{B_r}\(fv^++gv^-\)+\frac2r\int_{B_r}\[\(x\cdot\nabla
f\)v^++\(x\cdot\nabla g\)v^-\]\end{multline*} and hence
\begin{multline}\label{eq2.7}
\frac{\partial}{\partial r}\(\frac1{r^{n+2}}\int_{B_r}|\nabla
v|^2\)=-\frac4{r^{n+3}}\int_{\partial B_r}v\frac{\p v}{\p
\nu}\\+\frac{2(n+2)}{r^{n+3}}\int_{B_r}\(fv^++gv^-\)
+\frac2{r^{n+2}}\int_{\partial B_r}\(\frac{\p v}{\partial
\nu}\)^2-\frac2{r^{n+2}}\int_{\partial
B_r}\(fv^++gv^-\)\\+\frac2{r^{n+3}}\int_{B_r}\[\(x\cdot\nabla
f\)v^++\(x\cdot\nabla g\)v^-\],
\end{multline}
where we have also used the identity
\begin{multline*}
-\frac4{r^{n+3}}\int_{B_r}|\nabla v|^2=-\frac2{r^{n+3}}\int_{B_r}\[\Delta\(v^2\)-2v\Delta v\]
\\=-\frac4{r^{n+3}}\int_{\partial B_r}v\frac{\p v}{\p\nu}+\frac4{r^{n+3}}\int_{\partial B_r}\(fv^++gv^-\).
\end{multline*}
Since
\begin{multline*}\frac{\partial}{\partial
r}\(\frac2{r^{n+2}}\int_{B_r}\(fv^++gv^-\)\)\
\\=\frac{-2(n+2)}{r^{n+3}}\int_{B_r}\(fv^++gv^-\)+\frac2{r^{n+3}}\int_{\partial
B_r}\(fv^++gv^-\)\ d\sigma
\end{multline*}
and
\begin{equation*}
\frac{\partial}{\partial r}\(\frac2{r^{n+3}}\int_{\partial
B_r}v^2\)=-\frac8{r^{n+3}}\int_{\partial
B_r}v^2+\frac4{r^{n+3}}\int_{\partial B_r}v\frac{\p v}{\p\nu},
\end{equation*}
(\ref{eq2.5}) follows, with
$$e(r)=\frac2{r^{n+3}}\int_{B_r}\[\(x\cdot\nabla
f\)v^++\(x\cdot\nabla g\)v^-\].$$ (\ref{eq2.6}) is an immediate
consequence of this formula and the fact that $x_0\in S_v,\,v\in
C^{1,\gamma}$
\end{proof}

\begin{cor}\label{cor2.8}
If $f=f_0,\,g=g_0,$ both constants, and $W'(r)=0$ for $0<r<r_0,\,
v(x_0+x)$ is homogeneous of degree 2 in $x$.
\end{cor}

\begin{proof}
From the formula for $W'$ and the fact that $e\equiv0$ in this case.
\end{proof}

\begin{cor}\label{cor2.9}
$W_1(r)=W(r)+Dr^\gamma$ (where $D=D(n,\,\gamma,\,\|\nabla
f\|_\infty,\,\|\nabla
g\|_\infty,\,N)\geq0,\,D(-,\,-,\,0,\,0,\,-)\equiv0)$ is increasing
for $0<r<r_0.$
\end{cor}

For further use we will recall Kato's inequality:
\begin{lemma}[Kato \cite{K}]\label{lemma2.10} Assume that $w\in
W^{2,2}_{loc}(U)$. Then, $\Delta|w|\geq(\sign w)\Delta w$ in the
$H^1_{loc}(U)$ sense, i.e. for all $\theta\in
C_0^\infty(U),\,\theta\geq0$, we have
$$-\int\nabla|w|\cdot\nabla\theta\geq\int(\sign w)\Delta w\theta.$$
\end{lemma}

\begin{lemma}\label{lemma2.11}
For $0<r<r_0,\,x_0\in S_v,$ we have
$$\frac{\p}{\p r}\(\frac1{2r^{n+3}}\int_{\p B_r}v^2\)=\frac1r\[W_1(r)-\frac1{r^{n+3}}\int_{B_r}\[\(fv^++gv^-\)\]-Dr^\gamma\]$$
\end{lemma}

\begin{proof}
Recall from the proof of Lemma \ref{lemma2.4}, that
$$\frac{\p}{\p r}\(\frac1{2r^{n+3}}\int_{\p B_r}v^2\)=
-\frac2{r^{n+4}}\int_{\p B_r}v^2+\frac1{r^{n+3}}\int_{\p
B_r}v\frac{\p v}{\p \nu}.$$ But,
\begin{multline*} \int_{\p
B_r}v\frac{\p v}{\p \nu}=\frac12\int_{\p B_r}\frac{\p}{\p
r}\(v^2\)=\frac12\int_{B_r}\Delta\(v^2\)
\\=\int_{B_r}v\Delta v+|\nabla v|^2=\int_{B_r}|\nabla v|^2+\int_{B_r}\[fv^++gv^-\]
\end{multline*}
and the Lemma follows.
\end{proof}

We now let, for $0<r<r_0$, $\displaystyle
v_r(x)=\frac{v(rx+x_0)}{r^2},\,f_r(x)=f(rx+x_0),\,g_r(x)=g(rx+x_0)$,
where $x_0\in S_v.$ Note that $\Delta
v_r=f_r\chi_{\{v_r\geq0\}}-g_r\chi_{\{v_r<0\}}$ in
$B_1=B(0,r)\,(x_0=0).$

\begin{lemma}\label{lemma2.12}
Let
$v_r^{(1)}=f_rv_r^+(x)+\(g_r+\eta_0/2\)v^-(x)-a_1|x|^2,\,v_r^{(2)}=v_r^-+a_2|x|^2$,
where
$a_i=a_i(n,\,\tilde{B_1},\,\tilde{B_2},\,\tilde{B_3},\,N,\,\alpha,\,\Lambda)\geq0.$
Then:
\begin{enumerate}
\item[(i)] $v_r^{(1)}$ is superharmonic in $B_1.$
\item[(ii)] $v_r^{(2)}$ is subharmonic and non-negative in $B_1.$
\item[(iii)] $v_r^+$ is subharmonic in $B_1$.
\end{enumerate}
\end{lemma}

\begin{proof}
All functions are continuous, so we just need to check the sign of
the distributional Laplacian. Note that
$$v_r^{(1)}=\frac{\(f_r+g_r+\eta_0/2\)}2|v_r(x)|+\frac{\(f_r-g_r-\eta_0/2\)}2v_r(x)-a_1|x|^2$$
\begin{multline*}\Delta v_r^{(1)}=\frac{f_r+g_r+\eta_0/2}2\Delta\(|v_r(x)|\)+\frac{f_r-g_r-\eta_0/2}2\Delta v_r(x)\\+
2\frac{\nabla\(f_r+g_r\)}2\nabla\(|v_r|\)+2\frac{\Delta\(f_r-g_r\)}2\cdot
v_r\\+\frac{\Delta\(f_r+g_r\)}2|v_r|+ \frac{\Delta\(f_r-g_r\)}rv_r
-a_12n\\\leq \frac{\(f_r+g_r+\eta_0/2\)}2\(\sign
v_r\)\(f_r\chi_{\{v_r\geq0\}}-g_r\chi_{\{v_r<0\}}\)
\\+\frac{\(f_r-g_r-\eta_0/2\)}2\(\sign
v_r\)\(f_r\chi_{\{v_r\geq0\}}-g_r\chi_{\{v_r<0\}}\)+2\tilde{B_2}N
+2\tilde{B_2}N-a_12n\\
\leq\tilde{B_2}^2\tilde{B_2}^2+2\tilde{B_3}N-a_12n
\end{multline*}
and (i) follows. (Here we have used that $\(f_r+g_r+\eta_0/2\)<0$.)
Also,
$$v_r^{(2)}(x)=\frac{|v_r(x)|-v_r(x)}2+a_2|x|^2,$$
so that
\begin{multline*}
\Delta v_r^{(2)}(x)=\frac{(\sign
v_r)v_r(x)-\Delta v_r(x)}2+2na_2
\\=\frac{\(f_r\chi_{\{v_r\geq0\}}-g_r\chi_{\{v_r<0\}}\)-\(f_r\chi_{\{v_r\geq0\}}-g_r\chi_{\{v_r<0\}}\)}2+2na_2
\\\geq g_r\chi_{\{v_r<0\}}+2na_2
\end{multline*}
and (ii) follows. For (iii) note that $v_r^+=|v_r^+|+v_r^+/2$, so
that
\begin{multline*}
\Delta v_r^{+}(x)\geq\frac{(\sign v_r)\Delta v_r(x)+\Delta
v_r(x)}2+2na_2
\\=\frac{\(f_r\chi_{\{v_r\geq0\}}+g_r\chi_{\{v_r<0\}}\)+\(f_r\chi_{\{v_r\geq0\}}-g_r\chi_{\{v_r<0\}}\)}2
\\= f_r\chi_{\{v_r\geq0\}}\geq0.
\end{multline*}
\end{proof}

\begin{cor}\label{cor2.13}
$-\int_{B_1}\[f_rv_r^++\(g_r+\eta_0/2\)v_r^-\]\geq-a_3,$ where
$a_3>0$ has the same dependance as $a_i$ in Lemma \ref{lemma2.12}.
\end{cor}

\begin{proof}
$v_r^{(1)}$ is superharmonic in $B_1$, $v_r^{(1)}(0)=0$. Then
$$v_r^{(1)}(0)\geq\int_{B_1}f_rv_r^++\(g_r+\eta_0/2\)v_r^--a_1|x|^2\hspace{-6.02cm}-$$
and the Corollary follows.
\end{proof}

We now define, for $x_0\in S_v,\, 0<r<r_0,\,S(r)=\(\int_{\p
B_r}v^2\)^{1/2}.\hspace{-2.14cm}-$

\begin{lemma}[Non-degeneracy]\label{lemma2.14}
$\displaystyle\liminf_{r\to0}\frac{S(r)}{r^2}>0$
\end{lemma}

\begin{proof}
Assume without loss of generality, that $x_0=0.$ If the conclusion
fails, we can find $r_i\to0$ such that $\frac{S(r_i)}{r^2_i}\to 0.$
Let $v_i(x)=\frac{v(r_ix)}{r_i^2}$, so that $\int_{\p
B_1}v_i^2\to0.$ Note that, in $B_1$, $\Delta
v_i=f_{r_i}\chi_{\{v_i\geq0\}}-g_{r_i}\chi_{\{v_i<0\}}\geq\eta_0>0.$
Also, $|\Delta v_i|\leq2\tilde{B_1},$ $v_i(0)=0$. By subharmonicity
of $v_i$, $v_i(0)=0$ we see that
$\int_{B_1}v_i^-\leq\int_{B_1}v_i^+$. Since, by (iii) in Lemma
\ref{lemma2.12} $v_i^+$ is subharmonic, $\int_{B_1}v_i^+\leq
c_n\int_{\p B_1}v_i^+\leq c_n\(\int_{\p B_1}\(v_i^+\)^2\)^{1/2}.$
Thus,
$$\int_{B_1}|v_i|\leq\int_{B_1}v_i^++\int_{B_1}v_i^-\leq2\int_{B_1}v_i^+\leq 2c_n\(\int_{\p
B_1}\(v_i^+\)^2\)^{1/2}\to0.$$ After passing to a subsequence
$v_i\to v_0$, where the convergence is uniform on compact subsets of
$B_1$ and in $W_{loc}^{2,2}(B_1).$ But then $\Delta
v_0\geq\eta_0>0$, but $\int_{B_1}|v_0|=0$, a contradiction.
\end{proof}

\begin{rem}\label{rem2.15}
Note that the above proof shows that, if $S^+(r)=\(\int_{\p
B_r}\(v^+\)^2\)^{1/2},\hspace{-2.78cm}-\hspace{2.8cm}$ then
$\displaystyle\liminf_{r\to0}\frac{S^+(r)}{r^2}>0.$
\end{rem}

We now turn to the classification of blow-up points, following the
ideas of Monneau-Weiss \cite{MW}.

\begin{lemma}\label{lemma2.16}
Let $\displaystyle-M=\lim_{r\downarrow0}W_1(r).$ Assume that $x_0\in
S_{v}$ is such that $M<\infty$. Then, there exists
$G=G\(n,\,\tilde{B_1},\,\tilde{B_2},\,\tilde{B_3},N,\,M\)$ such that
$\displaystyle\sup_{0<r<r_0}\frac{S(r)}{r^2}\leq G.$
\end{lemma}

\begin{proof}
Note that, in view of Lemma \ref{lemma2.11}, if $0<r<r_0$ is such
that
\newline$\displaystyle-\frac1{r^{n+2}}\int_{B_r}\[fv^++gv^-\]>M+Dr^\gamma$,
then $\displaystyle\frac{\p}{\p r}\(\frac1{r^{n+2}}\int_{\p
B_r}v^2\)>0$. Note that the last inequality is equivalent with
$\displaystyle\frac{\p}{\p r}\(\int_{B_1}v_r^2\)>0$. Our first step
in the proof is to show that there exists
$C_1=C_1\(n,\,\tilde{B_1},\,\tilde{B_2},\,\tilde{B_3},N,\,\Omega\)$
such that, for $0<r<r_0$ we have

\begin{equation}\label{eq2.17}
\int_{\p B_1}\(v_r^+\)^2\leq C_1\left\{1+\int_{\partial
B_1}\(v_r^-\)^2\right\}
\end{equation}
In order to establish (\ref{eq2.17}), we first prove an auxiliary
Claim:

\begin{claim}\label{claim2.18}
For each $R>0$, there exists $\epsilon_0=\epsilon_0\(R,\,n\)$ such
that if $w(0)=0,\,\Delta w^+\geq0,\,0\leq\Delta
w\leq\epsilon_0,\,\int_{B_1}|\nabla w|^2\leq R$ and $\int_{\p
B_1}\(w^-\)^2\leq\epsilon_0,\,\(\int_{\p B_1}w^2\)^{1/2}\leq2$, then
$\int_{\p B_1}w^2\leq1/2$.
\end{claim}

\begin{proof}[Proof of Claim \ref{claim2.18}]
If not, we can find $R>0$ and functions $w_j$ with
$w_j(0)=0,\,0\leq\Delta w_j\leq1/j,\,\int_{B_1}|\nabla w_j|^2\leq
R,\,\(\int_{\p B_1}w_j^2\)^{1/2}\leq2,\,\int_{\p
B_1}\(w_j^-\)^2\leq1/j$ but $\int_{\p B_1}w_j^2\geq1/2.$ Since the
$w_j^+$ are subharmonic, $\int_{B_1}w_j^+\leq c_n\int_{\p
B_1}w_j^+\leq2c_n$. Since $w_j$ are subharmonic,
$w_j(0)=0,\,\int_{B_1}w_j^-\leq\int_{B_1}w_j^+\leq2c_n$. Hence, by
Poincare's inequality $\int_{B_1}w_j^2\leq\(R+4c_n\)\alpha_n.$
Hence, we can find a subsequence, (still called $j$) such that
$w_j\to w,$ uniformly on compact sets and $\int_{B_1}|\nabla
w|^2\leq R.$ Moreover, by compactness in the trace theorem, we have
$\int_{\p B_1}w^2\geq1/2$. We also have $\Delta
w=0,\,w(0)=0,\,\(\int_{\p B_1}w^2\)^{1/2}\leq2$ and $\int_{\p
B_1}w^-=0$. But then $w\geq0,\,w(0)=0$ and $\Delta w=0$ imply
$w\equiv0$, a contradiction.
\end{proof}

Suppose now that (\ref{eq2.17}) fails for some fixed $C_1>1$, to be
determined. Then, there exists a sequence $\{r_m\},\,0<r_m<r_0$ so
that $\int_{\p B_1}\(v_{r_m}^+\)^2\geq C_1\left\{1+\int_{\p
B_1}\(v_{r_m}^-\)^2\right\}.$ Using Corollary \ref{cor2.9}, we see
that
\begin{eqnarray}\label{eq2.19}
\nonumber\int_{B_1}|\nabla v_{r_n}|^2-2\int_{\p B_1}v_{r_n}^2&
\leq &
W_1(r_0)+2Dr_0^\gamma-2\int_{B_1}\(fv_{r_n}^++gv_{r_n}^-\)\\
& \leq & W_1(r_0)+2Dr_0^\gamma-2\int_{B_1}gv_{r_n}^-.
\end{eqnarray}

Consider now $w_n=v_{r_n}/\(\int_{\p B_1}\(v_{r_n}^+\)^2\)^{1/2}$.
Note that $w_n(0)=0,\,\Delta w_n\geq 0$ and by (iii) in Lemma
\ref{lemma2.12}, we have $\Delta w_n^+\geq0$. Also $\(\int_{\p
B_1}\(w_n\)^2\)^{1/2}\leq\(1+1/{C_1^{1/2}}\)\leq2,$ $\int_{\p
B_1}w_n^2\geq\int_{\p B_1}\(w_n^+\)^2=1$ and $|\Delta w_n|\leq
C/C_1^{1/2}$, where $C=C(\tilde{B_1})$, since $\int_{\p
B_1}\(v_{r_n}^+\)^2\geq C_1.$ Moreover, $\(\int_{\p
B_1}\(w_n^-\)^2\)^{1/2}\leq1/C_1^{1/2}.$

But, (\ref{eq2.19}) shows that
\begin{multline*}
\int_{B_1}|\nabla w_n|^2\leq2\int_{\p
B_1}w_n^2+\frac{W_1(r_0)}{\int_{\p
B_1}\(v_{r_n}^+\)^2}+\frac{2Dr_0^\gamma}{\int_{\p
B_1}\(v_{r_n}^+\)^2}\\+2c_n\tilde{B_1}\cdot\(\int_{\p
B_1}\(v_{r_n}^-\)+a_2\)\bigg/\int_{\p B_1}\(v_{r_n}^+\)^2,
\end{multline*} in view of Lemma \ref{lemma2.12}, (ii). Finally,
since $\int_{\p B_1}\(v_{r_n}^+\)^2\geq C_1\geq1,\, \(\int_{\p
B_1}w_n^2\)^{1/2}\leq2$ and $\int_{\p
B_1}\(v_{r_n}^-\)^2\leq\(\int_{\p
B_1}\(v_{r_n}^-\)^2\)^{1/2}\leq\frac1{C_1^{1/2}}\(\int_{\p
B_1}\(v_{r_n}^+\)^2\)^{1/2}$, we see that $\int_{B_1}|\nabla
w_n|^2\leq R,\,
R=R(\tilde{B_1},\,\tilde{B_2},\,\tilde{B_3},\,n,\,N,\,r_0)$, for all
$C_1\geq1.$ But if we now choose
$C/C_1\leq\epsilon_0,\,\frac1{C_1^{1/2}}\leq\epsilon_0,$ where
$\epsilon_0$ is as in Claim \ref{claim2.18}, we reach a
contradiction to Claim \ref{claim2.18}, establishing (\ref{eq2.17}).

We now proceed to the completion of the proof of Lemma
\ref{lemma2.16}. For $0<r<r_0,\,\tilde{r}\in(r/2,\,r),$ we have:
$W_1(r)-W_1(\tilde{r})\leq W_1(r_0)+M$. But, by Lemma
\ref{lemma2.4},
\begin{eqnarray*}
W_1(r)-W_1(\tilde{r})&=&\int_{\tilde r}^rW_1'(s)\
ds\\&=&\int_{\tilde r}^r2s \int_{\p B_1}\(\p_sv_s\)^2\ ds+
\int_{\tilde r}^re(s)\ ds+\gamma D\int_{\tilde r}^rs^{\gamma-1}\
ds\\&\geq&\int_{\tilde r}^r2s \int_{\p B_1}\(\p_sv_s\)^2\ ds,
\end{eqnarray*}
by our choice of $D$ and the fact that
$\displaystyle\p_sv_s=\frac{x\cdot\nabla
v(sx+x_0)}{s^3}-\frac{2v(sx+x_0)}{s^3}.$ The right hand side is
bigger than $\displaystyle r\int_{\tilde r}^r \int_{\p
B_1}\(\p_sv_s\)^2\ d\sigma ds$, which by Cauchy-Schwarz is bigger
than $\displaystyle\int_{\partial B_1}\(v_r-v_{\tilde r}\)^2\
d\sigma.$ Hence, for $0<r<r_0,\,\tilde r\in(r/2,\,r),$ we have
\begin{equation}\label{eq2.19'}
\int_{\partial B_1}\(v_r-v_{\tilde r}\)^2\leq\(W_1(r_0)+M\).
\end{equation}

We next show:

\begin{claim}\label{claim2.20}
There exists $\t M=\t
M(n,\,\tilde{B_1},\,\tilde{B_2},\,\tilde{B_3},\,N,\,M,\,r_0,\,\eta_0)$
such that if $\int_{\p B_1}v_r^2>\t M$, then $\frac{\p}{\p
r}\(\int_{\p B_1}v_r^2\)>0$, for $0<r<r_0.$
\end{claim}
To establish the Claim note that in light of the Remark at the
beginning of the proof of Lemma \ref{lemma2.16}, we only need to
show that $\int_{\p B_1}v_r^2\geq\t M$ implies that
$-\int_{B_1}\[f_rv_r^++g_rv_r^-\]>M+Dr^\gamma.$ By Corollary
\ref{cor2.13},
$-\int_{B_1}\[f_rv_r^++\(g_r+\eta_0/2\)v_r^-\]\geq-a_3$, so it
is enough to show that
\begin{equation}\label{eq2.21}
\eta_0\int_{B_1}v_r^->M+Dr^\gamma+a_3.
\end{equation}
From (ii) in Lemma \ref{lemma2.12}, we have (by interior estimates),
$$\(\int\limits_{1/2<|x|<3/4}\(v_r^-+a_2|x|^2\)^2\)^{1/2}\leq c_n\int_{B_1}\(v_r^-+a_2|x|^2\),$$
so that
\begin{equation}\label{eq2.22}
\int_{B_1}v_r^-\geq\frac1{c_n}\(\int\limits_{1/2<|x|<3/4}\(v_r^-\)^2\)^{1/2}-\tilde{c_n}.
\end{equation}
But, $\displaystyle\int\limits_{1/2<|x|<3/4}\(v_r^-\)^2\geq
a_n\int_{1/2}^{3/4}\int_{\p B_1}\(v_{rs}^-\)^2\ d\sigma ds$ and
$$\int_{\p B_1}\(v_{rs}^-\)^2=\int_{\p B_1}\(v_{rs}\)^2-\(v_{rs}^+\)^2
\geq\int_{\p B_1}\(v_{rs}\)^2-C_1-C_1\int_{\p B_1}\(v_{rs}^-\)^2,$$
from (\ref{eq2.17}). Thus, $\displaystyle\int_{\p
B_1}\(v_{rs}^-\)^2\geq\frac1{1+C_1}\int_{\p
B_1}\(v_{rs}\)^2-\frac{C_1}{C_1+1},$ and so, from (\ref{eq2.22}) we
obtain
$$\int_{B_1}v_r^-\geq d_n\(\int_{1/2}^{3/4}\int_{\p B_1}v_{rs}^2\ d\sigma ds\)-C_2,$$
with $C_2$ having the same dependence as $C_1$. If we now use
(\ref{eq2.19}) with $\t r=rs,$ we see using (3.9),
$$\int_{B_1}v_r^-\geq\t{d_n}\(\int_{\p B_1}v_{r}^2\)^2-b_n\(W_1(r_0)+M\)^{1/2}-C_2$$
and (\ref{eq2.21}) holds for $\t M$ large enough.

We can now conclude the proof of Lemma \ref{lemma2.16}: if
$S(r)/r^2\leq\t M$ for $0<r<r_0$, we are done. If for all
$0<r<r_0,\,S(r)/r^2>\t M$, then by Claim \ref{claim2.20} we have
$S(r)/r^2=\int_{\p B_1}v_r^2<S(r_0)/r_0^2$ for $0<r<r_0$ and we are
also done. Note that, if for some $0<r_1<r_0$ we have
$S(r_1)/r_1^2>\t M,$ then for all $r_1<r<r_0$ we have $S(r)/r^2>\t
M,$ by virtue of Claim \ref{claim2.20}. It is now easy to show that
$S(r)/r^2\leq\max\(\t M,\,S(r_0)/r_0\),\,0<r<r_0.$ Thus, Lemma
\ref{lemma2.16} follows.

\end{proof}

\begin{cor}\label{cor2.23}
Let $M,\,G$ be as in Lemma \ref{lemma2.16}. Then there exists $\t
G$, with the same dependance as $G$ such that, for $0<r<r_0/2$, we
have
$$\sup_{|x|\leq1}|v_r(x)|+\(\int_{B_1}|\nabla v_r|^2\)^{1/2}\leq\t G.$$
\end{cor}

\begin{proof}
By Corollary \ref{cor2.9}, for $0<r<r_0$ we have
\begin{eqnarray*}
\int_{B_1}|\nabla v_r|^2&\leq&2\int_{\p
B_1}v_r^2-2\int_{B_1}\(fv_r^++gv_r^-\)+2Dr_0^\gamma+W_1(r_0)\\
&\leq&2\int_{\p
B_1}v_r^2-2\int_{B_1}gv_r^-+2Dr_0^\gamma+W_1(r_0).
\end{eqnarray*}
Using now Lemma \ref{lemma2.12} (ii) and Lemma \ref{lemma2.16}, the
gradient estimate follows. For the $L^\infty$ estimate, we use Lemma
\ref{lemma2.12} (i) and (ii) and the fact that for non-negative
subharmonic functions, the $L^2$ spherical averages are increasing.
Thus, for instance

\begin{multline*} \sup_{|x|\leq1}|v_r^+(x)|\leq\sup_{|x|=1}|v_r^+(x)|
\leq \t{c_n}\(\int_{\p B_1}|v_r^+|^2\)^{1/2}\\\leq
c_n\(\int\limits_{1/2<|x|<3/2}|v_r^+|^2\)^{1/2}\leq\t c_n\(\int_{\p
B_1}\(v_{2r}^+\)^2\)^{1/2}\end{multline*} and correspondingly for
$v_r^-$.
\end{proof}

We are now ready, in analogy with \cite{MW}, to state our
classification of blow-up points.

\begin{thm}\label{thm2.24}
Assume that $x_0\in S_v$ and $W_1(r)$ is defined in Corollary
\ref{cor2.9}.
\begin{enumerate}
\item[(i)] If
$\displaystyle\lim_{r\downarrow0}W_1(r)=-M,\,M<+\infty$, then
$S(r)/r^2$ and $\|v_r\|_{W^{2,p}(B_1)},\,1<p<\infty$ remain bounded
for $0<r<r_0/2.$ Moreover, if $\{r_j\}$ is a sequence tending to 0,
after passing to a subsequence $\{r_{j'}\}$, the functions
$v_{r_{j'}}$ converge in $C^{1,\gamma}(B_1),\,0\leq\gamma<1$ and
$W^{2,p}(B_1),\,1\leq p<\infty$ to a function $\bar v.$ The function
$\bar v$ solves the equation
$$\Delta\bar v=f_0\chi_{\{\bar v\geq0\}}-g_0\chi_{\{\bar v<0\}}\text{ in }\R^n,$$
with $f_0=f(x_0),\,g_0=g(x_0)$ and is homogeneous of degree 2.

\item[(ii)] If $\displaystyle\lim_{r\downarrow0}W_1(r)=-\infty$,
then $\displaystyle\lim_{r\downarrow0}\frac{S(r)}{r^2}=+\infty$. Let
$r_j\downarrow0$ and define $\displaystyle
w_j(x)=\frac{v(r_jx+x_0)}{S(r_j)},\,T_j=\frac{S(r_j)}{r_j^2}.$ Then,
after passing to a subsequence $\{r_{j'}\}$, $w_j$ converge in
$C^{1,\gamma}(B_1)$ and $W^{2,p}(B_1),\,0\leq\gamma<1,\,1\leq
p<\infty$ to a harmonic function $\bar w,$ with $\bar
w(0)=\nabla\bar w(0)=0$, which is non-zero and homogeneous of degree
2.
\end{enumerate}
\end{thm}

\begin{proof}
From Corollary \ref{cor2.23}, in case (i) it only remains to show
that $\bar v$ is homogeneous of degree 2. But, for any $0<s<1$ we
have $W(s)=W(s;\,0;\,\bar v)=\displaystyle\lim_{j'\to\infty}
W_1(sr_{j'},\,x_0,\,v)=-M.$ Then, from Corollary \ref{cor2.8} $\bar
v$ is homogeneous of degree 2.

For case (ii), we must have
$$\lim_{r\downarrow0}\int_{B_1}|\nabla v_r|^2+2\int_{B_1}fv_r^2+2\int_{B_1}gv_r^--2\int_{\p B_1}v_r^2=-\infty$$
But then, since $f>0,\,g<0$, we must have that

\begin{equation}\label{eq2.25}
\lim_{r\to0} 2\int_{\p B_1}v_r^2-2\int_{B_1}gv_r^-=+\infty.
\end{equation}
By Lemma \ref{lemma2.12} (ii),
$$\int_{B_1}v_r^-\leq
c_na_2+c_n\int_{\p B_1}v_r^-\leq c_na_2+\(\int_{\p
B_1}\(v_r^-\)^2\)^{1/2}.$$ But then, since $(-g)\geq\eta_0,$ and
$-g\leq\t B_1$, we conclude from (\ref{eq2.25}) that
$\displaystyle\lim_{r\to0}\int_{\p
B_1}v_r^2+\int_{B_1}v_r^-=+\infty$, which in turn implies
$\displaystyle\lim_{r\to0}\int_{\p B_1}v_r^2=+\infty$, or
$\displaystyle\lim_{r\to0}\frac{S(r)}{r^2}=+\infty.$ By Corollary
\ref{cor2.9}, dividing by $T_j^2$, we obtain

\begin{equation}\label{eq2.26}
\int_{B_1}|\nabla
w_j|^2\leq\frac{W_1(r_0)}{T_j^2}+\frac2{T_j}\int_{B_1}\[f_{r_j}w_j^++g_{r_j}w_j^-\]+2\int_{\p
B_1}w_j^2-\frac{Dr_j^\gamma}{T_j^2}.
\end{equation}
Also, for $j$ large, $|\Delta w_j|\leq 1$ in $B_1,\,\int_{\p
B_1}w_j^2=1,\,\Delta w_j^+\geq0,\,\Delta w_j\geq0,\,w_j(0)=0.$ Then,
$\int_{B_1}(w_j)^2\leq C$ and from the formulae above
$\int_{B_1}|\nabla w_j|^2\leq3$, for $j$ large. Thus, the $w_j,$
after passing to a subsequence $j',$ converge uniformly on compacts
and in $C^{1,\gamma}(B_1),\,W^{2,p}(B_1)$ to a $\bar w$ which is
harmonic in $B_1,\,\bar w(0)=\nabla\bar w(0)=0.$ Also, by
compactness of the trace operator, $\int_{\p B_1}\bar w^2=1$, so
that $\bar w$ is not zero. But, from (\ref{eq2.26}), we conclude
that $\int_{B_1}|\nabla\bar w|^2\leq2\int_{\p B_1}|\bar w|^2$. But
then, by the Almgren monotonicity formula (see for example Lemma 4.2
in \cite{MW}) $w$ is homogeneous of degree 2.
\end{proof}

\begin{cor}[No mixed asymptotics]\label{cor2.27}
We cannot have for two sequences $\{r_j\},\,\{\t r_j\},$ both
tending to 0 that $\displaystyle
\lim_{j\to\infty}\frac{S(r_j)}{r_j^2}=+\infty$, but $\displaystyle
\sup_j\frac{S(\t r_j)}{\t r_j^2}<+\infty$.
\end{cor}

\begin{proof}
If $\displaystyle \lim_{r\downarrow0}W_1(r)=-\infty$, then for all
such sequences the limit is $+\infty$. On the other hand if
$\displaystyle \lim_{r\downarrow0}W_1(r)>-\infty$, we have
boundedness near $r=0$. In either case the mixed asymptotic
assumption leads to a contradiction.
\end{proof}

We will next use these results to study partial regularity of the
free boundary $\F.$ We start with a 2-dimensional result, due to
Shahgholian (\cite{S})

\begin{thm}[\cite{S}]\label{thm2.28}
Let $v$ be the solution of (\ref{eq2.2}), when $n=2$, under our
assumptions. Assume that $x_0\in S_v$ is such that $|\{v<0\}\cap
B(x_0,r)|\geq c_0r^2$, for $0<r<r(x_0)$, with $c_0>0.$ Then $x_0$ is
an isolated point of $S_v.$
\end{thm}
We will provide a proof of this Theorem, (following \cite{S}), for
the reader's convenience. The key point is the following

\begin{lemma}[\cite{S}]\label{lemma2.29}
Assume that $\bar v$ is a homogeneous of degree 2 solution to
(\ref{eq2.2}) in $\R^2$, with $f=f_0,\,g=g_0$, both constants. (As
before, $f_0>0,\,g_0<0,\,f_0+g_0<0$). Then, $S_{\bar v}=\{0\}$, or,
after rotation, $S_{\bar v}=\{(x_1,x_2)=(0, x_2):x_2\in\R\}.$ In
this case $\bar v=\frac{f_0}2x_1^2.$
\end{lemma}

\begin{proof}
Recall that $\Delta\bar v\geq\eta_0>0$ ($\eta_0=\min(f_0,\,g_0)$).
Assume that $S_v\neq\{0\}$. After rotation we can assume that, by
the homogeneity of $\bar v$, $(0,\,1)\in S_{\bar v}$, so that
$\lambda(0,\,1)\in S_{\bar v},\,\lambda>0$. Assume first that $\bar
v\geq0$ in a neighborhood of (0, 1). Then, in an angle, $\Delta\bar
v=f_0.$ Consider $w=\bar v-\frac{f_0}2x_1^2$. Then, in this angle,
by uniqueness in the Cauchy problem, $w\equiv0$. But this argument
can be continued all around, so that $\bar v=\frac{f_0}2x_1^2$.
Thus, if not, there exists a neighborhood of (0, 1) in which $\bar
v<0$ is non-empty. Assume, for instance that the negative point is
in the top right quadrant.  By homogeneity, the point can be taken
on the unit circle. But then, all the points in the unit circle
between this point and the vertical axis are points where $\bar v$
is negative, otherwise we would have a local maximum, contradicting
the subharmonicity of $\bar v$. But then, if we consider a small
half-ball in the top right quadrant, centered at (0, 1), the Hopf
maximum principle yields a contradiction to $\bar
v(0,\,1)=0,\,\nabla\bar v(0,\,1)=0.$
\end{proof}

\begin{proof}[Proof of Theorem \ref{thm2.28}]
We can assume, without loss of generality, that $x_0=0$. Suppose we
have $x_j\in S_v,\,x_j\to0$. Let $r_j=|x_j|.$ Assume first that
$\ds\lim_{r\downarrow0}W_1(r)=-\infty.$ Then, by Theorem
\ref{thm2.24} (ii), $\frac{v(r_jx)}{S(r_j)}$, after passing to a
subsequence, converges in $C^{1,\gamma}(B_1),\,L^2(\p B_1)$ to a
harmonic polynomial $\bar w$ homogeneous of degree 2 and non-zero.
Moreover, $\frac{x_j}{|x_j|}\to\bar x\in\p B_1$, and $\bar w(\bar
x)=0,\,\nabla\bar w(\bar x)=0.$ But, when $n=2$, $\bar w$ must be a
rotate of $a\(x_1^2-x_2^2\)$ and hence $S_{\bar w}=\{0\}$, a
contradiction. If $\ds\lim_{r\downarrow0}W_1(r)>-\infty$, by Theorem
\ref{thm2.28} (i), $\frac{v(r_jx)}{r_j^2}$ converges, after passing
to a subsequence to a $\bar v$, a homogeneous of degree 2 solution,
$f=f_0,\,g=g_0$. Clearly $|\{\bar v<0\}\cap B_1|\geq c_0$. Also,
$\bar x\in S_{\bar v}$, so that by Lemma \ref{lemma2.29}, $\bar
v=\frac{f_0}2x_1^2$, after a rotation, which is a contradiction.
\end{proof}

We will extend next Theorem \ref{thm2.28} to $n>2$. The argument is
a standard one from the theory of minimal surfaces (see Chapter 11
of \cite{G}, whose notation for Hausdorff measures and Hausdorff
dimension we adopt). Similar arguments have been  used by Weiss
\cite{W} and Monneau-Weiss \cite{MW} in the context of free boundary
problems. Our result here is:

\begin{thm}\label{thm2.30}
Let $v$ be a solution of (\ref{eq2.2}), $n\geq2$, under our
assumptions. Let $\t S_v=\left\{x_0\in S_v:|\{v<0\}\cap
B(x_0,r)|\geq c_0r^n, \text{ for }0<r<r_0(x_0)\right\}.$ Then, for
each fixed $c_0>0$, the Hausdorff dimension of $\tilde S_v$ is at
most $n-2.$
\end{thm}

\begin{proof}
Fix $k>n-2$, we need to show that $H_k\(\t S_v\)=0$. Assume not, so
that $H_k\(\t S_v\)>0$. Consider the sets
$$\t S_v^j=\left\{x_0\in
S_v:|\{v\leq0\}\cap B(x_0,\,r)|\geq c_0r^n, \text{ for
}0<r<1/j\right\}.$$ Then, $\ds \t S_v=\bigcup_{j=j_0}^\infty\t
S_v^j$, where $1/j_0<r_0.$ Then, for some $\bar j\geq j_0$ we have
$H_k\(\t S_v^{\bar j}\)>0$. Hence by Proposition 11.3 in \cite{G},
for $H_k$-almost all $x_0\in\t S_v^{\bar j}$, we have
\begin{equation}\label{eq2.31}
\limsup_{r\to0}\frac{H_k^\infty\(\t S_v^{\bar j}\cap B(x_0,\,
r)\)}{\omega_kr^k}\geq 2^{-k}.
\end{equation}
Fix such an $x_0$, which we assume, without loss of generality, to
be 0. Choose a sequence $r_n\to 0$ such that
$$\frac{H_k^\infty\(\t S_v^{\bar j}\cap B_{r_n}\)}{\omega_kr_n^k}\geq 2^{-k}.$$
Consider $v_n(x)=\frac{v(r_nx)}{S(r_n)}$ and let $\bar v(x)$ be a
blow-up limit of a subsequence of $v_n$, in the sense of Theorem
\ref{thm2.24}. Fix a compact set $K$ in $B_1,\, U$ open $\subset
B_1,$ with $U\supset K\cap\t S_{\bar v}^{\bar j}$. Assume that
$x_n\in \t S_{v_n}^{\bar j},\,x_n\in K\setminus U$ and after passing
to a subsequence, assume that $x_n\to\bar x\in K\setminus U.$ Then,
$v_n(x_n)\to\bar v(\bar x),\,\nabla v_n(x_n)\to \nabla\bar v(\bar
x)$, so that $\bar x\in S_{\bar v}.$ Also, fix $0<r<1/\bar j$. Then
$$|\{\bar v\leq0\}\cap B(\bar x,\,r)|=|\{\bar v<0\}\cap B(\bar x,\,r)|
=\lim_{n\to\infty}|\{v_r<0\}\cap B(x_n,\,r)|\geq c_0r^n,$$ and so
$\bar x\in\t S_{\bar v}^{\bar j}$, but $\bar x\in K\setminus U$, and
$K\cap S_{\bar v}^{\bar j}\subset U$, which is a contradiction.
Thus, we have shown that there exists $n_0$ so that, for $n>n_0$, we
have
\begin{equation}\label{eq2.32}
U\supset K\cap \t S_{v_n}^{\bar j}
\end{equation}
Then, the proof of Lemma 11.5 in \cite{G} shows that for all
$K\Subset B_1,$ we have
\begin{equation}\label{eq2.33}
H_k^\infty\(K\cap\t S_{\bar v}^{\bar
j}\)\geq\limsup_{n\to\infty}H_k^\infty\(K\cap\t S_{v_n}^{\bar j}\).
\end{equation}
We next claim that
\begin{equation}\label{eq2.34}
\left\{x/r_n:x\in\t S_v^{j}\right\}\subset\t S_{v_{r_n}}^{\t j}.
\end{equation}
In fact, clearly $v_n(x/r_n)=0,\,\nabla v_n(x/r_n)=0$. Consider now
$\{y:v_n(y)<0\}\cap B(x/r_n,\,r),\,0<r<1/\bar j.$ This equals
$\{y:v_n(y)<0\}\cap \{y:|y-x/r_n|<r\}$. By the transformation
$y=z/r_n$, this set equals
$$\{z:v(z)<0\}\cap \left\{z:\left|\frac z{r_n}-\frac x{r_n}\right|<r\right\}
=\{z:v(z)<0\}\cap \left\{z:\left|z-x\right|<rr_n\right\}.$$ Also, if
$0<r<1/\bar j,\, rr_n<1/\bar j,\, n$ large. The Lebesgue measure of
the set of $y$'s equals $(r_n)^{-n}$ times that Lebesgue measure of
the set of $z$'s, which is then bigger than $\ds\frac1{(r_n)^n}\cdot
c_0(rr_n)^n=c_0r^n$, so that $x/r_n\in\t S_{v_r}^{\bar j}$. But
then,
$$H_k^\infty\(B_1\cap\t S_{v_n}^{\bar j}\)\geq \frac{H_k^\infty\(B_{r_n}\cap\t S_{v}^{\bar j}\)}{\omega_kr_n^k}
\geq2^{-k,}$$ by our choice of $r_n$. Hence, using (\ref{eq2.33}),
we see that
\begin{equation}\label{eq2.34'}
H_k^\infty\(B_1\cap\t S_{\bar v}^{\bar j}\)>0
\end{equation}

We now consider our classification of blow-ups. If
$\ds\lim_{r\downarrow0}W_1(r)=-\infty$, then, by (ii) $\bar v$ is a
non-zero, homogeneous of degree 2 harmonic polynomial. But then, as
is well-know $H_{n-2}(S_{\bar v})<\infty,\,S_{\bar v}\supset\t
S_{\bar v}^{\bar j},$ which contradicts (\ref{eq2.34'}) since
$k>n-2.$ If $\ds\lim_{r\downarrow0}W_1(r)>-\infty$, in view of
Theorem \ref{thm2.24} (i) and Lemma \ref{lemma2.14}, after passing
to a further sequence, we can assume that
$\frac{r_n^2}{S(r_n)}\to\alpha,\,\alpha\in(0,\,\infty).$ Hence,
$\alpha\bar v=\bar v_1$, where $\bar v_1$ is a homogeneous of degree
2 solution to (\ref{eq2.2}) with $f=f_0,\,g=g_0$, both constants. We
can now do the dimension reduction. From (\ref{eq2.34'}), we know
that $H_k^\infty\(B_1\cap\t S_{\bar v}^{\bar j}\)>0$. Using Lemmas
11.2 and 11.3 in \cite{G}, we can find $\bar x\in\t S_{\bar v}^{\bar
j}\setminus\{0\}$ such that $\ds\lim_{r\to0}\frac{H_k^\infty\(\t
S_{\bar v}^{\bar j}\cap B(\bar x,\,r)\)}{\omega_kr^k}\geq2^{-k}.$ By
homogeneity of $\bar v_1$, we can assume that $\bar x\in\p B_1.$ We
can pick a sequence $r_n\to0$, and consider a blow-up limit $\bar
v_{1,0}$, at $\bar x$, with respect to $r_n$. By the homogeneity of
$\bar v_1$, it is easy to see that $\bar v_{1,0}$ is constant in the
$\bar x$ direction. After rotation, we can assume this direction to
be the $x_n$ direction. But, it is easy to see that
$(x_1,\,x_2,\cdots x_{n-1},\,x_n)\in\t S_{\bar
v_{1,0}|_{\R^{n-1}}}^{\bar j}$ and that $H_{k-1}\(\t S_{\bar
v_{1,0}|_{\R^{n-1}}}^{\bar j}\)>0$. Proceeding in this way $n-2$
times, we find a contradiction to Theorem \ref{thm2.28}, which
concludes the proof.
\end{proof}

We are now ready to establish partial $C^{1,1}$ bounds.

\begin{defi}\label{defi2.35}
Let $f$ be a $C^{1,\gamma},\,0\leq\gamma<1$ function defined in a
neighborhood of a point $x_0$. We say that $f$ satisfies $C^{1,1}$
bounds at $x_0$ if
$$\lim_{r\to0}\sup_{|x-x_0|\leq r}\frac{|f(x)-(x-x_0)\nabla f(x_0)-f(x_0)|}{r^2}<+\infty.$$
We call the above limit ``the $C^{1,1}$ norm of $f$ at $x_0$".
\end{defi}

Our next task is to show that our solutions $v$ verify $C^{1,1}$
bounds at all $x_0\in\F$, except for a set of Hausdorff dimension at
most $n-2$. We start out with some preliminary results.

\begin{lemma}\label{lemma2.36} There exists a constant $c_n$ such
that for all homogeneous of degree 2 harmonic polynomials
$p,\,p\not\equiv0,$ we have
$$|\{p<0\}\cap B_1|\geq c_n.$$
\end{lemma}

\begin{proof}
We can assume $\int_{B_1}p^2=1$. If the conclusion fails, we can
find a sequence $p_j,\,\int_{B_1}p_j^2=1,\, p_j$ a harmonic
polynomial, homogeneous of degree 2, with $|\{p_j<0\}\cap
B_1|\xrightarrow[j\to0]{}0$. After passing to a subsequence, $p_j\to
p_0,\,p_0$ a harmonic polynomial, homogeneous of degree 2,
$\int_{B_1}p_0=1$ and such that $|\{p_0<0\}\cap B_1|=0$. By
homogeneity, $p_0\geq0$, but $p_0(0)=0$, so that $p_0\equiv0$, a
contradiction.
\end{proof}

\begin{lemma}\label{lemma2.37}
Let $c_n$ be as in Lemma \ref{lemma2.36}. Assume that $v$ is a
solution, $x_0\in S_v.$ Assume that for some sequence $r_j\to0$,
$\ds\sup_{|x-x_0|<r_j}\frac{|v(x)|}{r_j^2}\to\infty.$ Then,
$$|\{v<0\}\cap B(x_0,r)|\geq\frac{c_n}2r^n,\text{ for }0<r<r_0(x_0).$$
\end{lemma}

\begin{proof}
If not, there exists $\t r_j\to0$, such that
$$|\{v<0\}\cap B(x_0,\t r_j)|<\frac{c_n}2\(\t r_j\)^n.$$
But, by the proof in Corollary \ref{cor2.23}, we see that
$\frac{S(2\t r_j)}{(2\t r_j)^2}\to+\infty$. By Corollary
\ref{cor2.27} we have $\frac{S(\t r_j)}{\t r_j^2}\to+\infty$. But
then, by Theorem \ref{thm2.24} (ii), $\frac{v(\t r_j x+x_0)}{S(\t
r_j)}$ converges, after passing to a subsequence, to a $\bar w$
which is a non-zero harmonic polynomial homogeneous of degree 2. But
then, $|\{\bar w<0\}\cap B_1|\leq\frac{c_n}2$, which contradicts
Lemma \ref{lemma2.36}.
\end{proof}

\begin{thm}[Pointwise $C^{1,1}$ bounds on $S_v$]\label{thm2.38} Let
$v$ be a solution. Consider the set $B_v=\{x_0\in S_v:v\text{ does
not have pointwise $C^{1,1}$ bounds at $x_0$}\}$. Then, the
Hausdorff dimension of $B_v$ is at most $n-2.$
\end{thm}

\begin{proof}
Combine Lemma \ref{lemma2.37} with Theorem \ref{thm2.30}.
\end{proof}

\begin{rem}\label{rem2.39}
If $x_0\in\F,\,\nabla v(x_0)\not\equiv0$, then by \cite{CGK} $\F$ is
real analytic in a neighborhood of $x_0$ and by boundary elliptic
regularity we obtain $C^{1,1}$ bounds at $x_0$, Thus, the set of
points in $\F$ for which $v$ does have pointwise $C^{1,1}$ bounds
has Hausdorff dimension at most $(n-2).$
\end{rem}

\begin{rem}\label{rem2.40}
The results in Theorems \ref{thm2.28}, \ref{thm2.30}, \ref{thm2.38}
and in Remark \ref{rem2.39} are sharp. We show this for the case
$f=f_0,\,g=g_0$ constants. In order to show this, we make some
preliminary comments, in the case $n=2.$ In this case, Blank
(\cite{B}) found all homogeneous of degree 2 solutions, for which
$\{v<0\}\neq\emptyset.$ The calculation in Appendix 2 shows that,
for these solutions, $W(1)>-A$, where $A$ depends only on
$f_0,\,g_0$. Shahgholian (\cite{S}) observed that there are other
homogeneous of degree 2 solutions, which are non-negative. In fact,
any such solution $\bar v$ verifies $\Delta\bar v=f_0,\,\bar v\geq0$
in $R^2.$ Let $w=\bar v-\frac{f_0}4(x_1^2+x_2^2).$ This is a
harmonic polynomial, homogeneous of degree 2, so that, after
rotation $w=a(x_1^2-x_2^2)$ or $\bar
v=\(a+\frac{f_0}4\)x_1^2+\(\frac{f_0}4-a\)x_2^2$. Since $\bar
v\geq0$, we must have $-\frac{f_0}4\leq a\leq\frac{f_0}4.$ For those
solutions we also find $W(1)>-A,\,A$ depending only on $f_0,\,g_0$.
Combining these comments with Theorem \ref{thm2.24}, we se that, for
$n=2$ there exists $A=A(f_0,\,g_0)$ such that, if for $v$ we have
$\ds\lim_{r\downarrow0}W_1(r)<-A$, then
$\ds\lim_{r\downarrow0}W_1(r)=-\infty$ and
$\ds\lim_{r\downarrow0}\frac{S(r)}{r^2}=+\infty.$ One can then use
the argument in \cite{AW} to see that, by the Andersson-Weiss
construction we can find solutions (taking $M$ large in \cite{AW})
so that $W_1(1)<-A$, and hence, solutions which don't have $C^{1,1}$
bounds towards 0. In light of Lemma \ref{lemma2.36}, this shows the
sharpness of Theorem \ref{thm2.28} and of Theorem \ref{thm2.38} when
$n=2.$ To create higher dimensional examples, one just adds $n-2$
dummy variables. It remains a challenging problem to see if such
pathology can hold for solutions of (\ref{eq2.2}).
\end{rem}

We now turn to the issue of uniform pointwise $C^{1,1}$ bounds.

\begin{thm}\label{thm2.41}
Let $S_v^{(1)}=S_v/S_v^{(2)}$, where $$S_v^{(2),j}=\{x_0\in
S_v:|\{v<0\}\cap
B(x),r)|\geq\frac1jr^n,\,0<r<r_{0,j}(x_0)\},\,S_v^{(2)}=\bigcup_{j=1}^\infty
S_v^{(2),j}.$$

Note that Theorem \ref{thm2.30} shows that the Hausdorff dimension
of $S_v^{(2)}$ is at most $n-2$. Then, for $x_0\in S_v^{(1)}$ we
have uniform $C^{1,1}$ estimates, i.e. there exists $C=C(\t B_1,\,\t
B_2,\,\t B_3,\,n,\,\eta_0,\,r_0,\,N)>0$ such that for all $x_0\in
S_v^{(1)},\ds \sup_{\substack{ |x-x_0|\leq
r\\0<r<r_0/2}}\frac{|v(x)|}{r^2}\leq C.$
\end{thm}

\begin{proof}
In light of Theorem \ref{thm2.24}, Lemma \ref{lemma2.36} and
Corollary \ref{cor2.23} it suffices to show that for such $x_0$
$\ds\lim_{r\downarrow0}W_1(r)>-A$, where $A$ has the right
dependence. Let $\bar v$ be a blow-up limit at such an $x_0.$
Clearly, $\bar v\geq0.$ Thus, it suffices to show that, for such
$\bar v,\,W(1,\,\bar v)>-A.$ But, $\Delta\bar
v=f_0,\,\int_{B_1}\Delta\bar v=w_nf_0=\int_{\p B_1}\frac{\p\bar
v}{\p\nu}=2\int_{\p B_1}\bar v$, since $\bar v$ is homogeneous of
degree 2. Thus, $\int_{\p B_1}\bar v=\frac{w_nf_0}2$. Since $\bar v$
is non-negative and subharmonic, $\int_{B_1}\bar v\leq c_nf_0w_n/2.$
The rest of the proof follows easily from interior estimates and
homogeneity.
\end{proof}

\begin{rem}\label{rem2.42} Similarly, if $K\Subset\{x_0\in\F:\nabla v(x_0)\neq0\}$
we also have uniform pointwise $C^{1,1}$ bounds on $K$. (See Remark
\ref{rem2.39}).
\end{rem}

Our final result is a partial regularity result for $\F.$

\begin{thm}\label{thm2.43}
Let $v$ be a solution of (\ref{eq2.2}) satisfying our assumptions.
Then $\F=\F_0\cup S_v^{(1)}\cup S_v^{(2)}$, where $S_v^{(2)}$ has
Hausdorff dimension at most $(n-2)$, $S_v^{(1)}$ is $(n-1)$ regular
i.e. $H_{(n-1)}\(S_v^{(1)}\)\leq C$, with $C=C(\t B_1,\,\t B_2,\,\t
B_3,\,N,\,\eta_0,\,r_0,\,n)$ and $\F_0$ is relatively open and for
each $x_0\in\F$ there exists a neighborhood $U_{x_0}$ such that
$\F\cap U_{x_0}$ is a real-analytic hypersurface.
\end{thm}

\begin{proof}
$\F_0=\{x_0\in\nabla v(x_0)\neq0\}$ and $S_v^{(1)},\,S_v^{(2)}$ are
defined in Theorem \ref{thm2.41}. From Theorem \ref{thm2.41} we know
that the Hausdorff dimension of $S_v^{(2)}$ is at most $n-2$, so it
remains to show that $S_v^{(1)}$ is $(n-1)$ regular, (in light of
Theorem 8 in \cite{CGK}, which shows the desired property of
$\F_0$). In order to show this, we make some preliminary claims.

\begin{claim}\label{claim2.44}
If $x_0\in S_v^{(1)}$ (without loss of generality, we take $x_0=0$)
we have, for $0<r<r_0/4,\,x\in B_r$, $|\nabla v(x)|\leq Cr$, with
$C$ as in the statement of Theorem \ref{thm2.43}.
\end{claim}

In order to establish the Claim, note that for $x\in B_{2r}$, we
have $|v(x)|\leq C|x|^2$, by Theorem \ref{thm2.41}. Next, we use
Lemma \ref{lemma2.12} (ii) and (iii) to obtain:
$$\int_{B_r}|\nabla v^+|^2\leq c_nCr^{n+2}\text{ and }
\int_{B_r}|\nabla v^-|^2\leq c_n\{C+a_2\}r^{n+2}$$ so that
$\int_{B_r}|\nabla v|^2\leq c_n(C+a_2)r^{n+2}$. Next, consider
$v_r(x)$ on $B_1$. We have $\int_{B_1}|v_r|^2\leq
C,\,\int_{B_1}|\nabla v_r|^2\leq C,\,$ and $|\Delta v_r|\leq C.$
 From this it is easy to see that, for $|x|\leq1/2$ we have $|\nabla v_r|\leq
 C$, which is our claim.
 \end{proof}

The next step is:
\begin{claim}\label{claim2.45}
Let $x_0\in S_v,\,e_i$ be a fixed coordinate direction,
$v_{e_i}=e_i\cdot\nabla v.$ Then, for $0\leq h$, small, we have
$$\int\limits_{B(x_0,r_0/2)\cap\{x:|\nabla v|\leq h\}}|\nabla v_{e_i}|^2\leq Ch.$$
\end{claim}

To establish Claim \ref{claim2.45}, we first introduce a truncation
of $v_{e_i}\in W^{1,2}(U)\cap C^{\gamma}\(\overline U\)$, to obtain
$\bar v_{e_i},$ where
$$\bar v_{e_i}=\begin{cases}
v_{e_i} & \text{if $-h<v_{e_i}<-\delta$ or $\delta<v_{e_i}<h,$}\\
0 & \text{if }|v_{e_i}|\leq\delta,\\
h & \text{if }|v_{e_i}|\geq h.
\end{cases}$$

Let $\psi$ be a standard mollifier and for $0<\epsilon\ll \delta$,
consider the mollifier $v_{e_i}*\psi_\epsilon$. We will apply
Green's Theorem to
$$\int_{B_r}\nabla\bar v_{e_i}\cdot\nabla\(v_{e_i}*\psi_\epsilon\),\text{ for } \frac{r_0}2<r<r_0,$$
where we have assumed, without loss of generality that $x_0=0$.
Since $|\F|=0$ (see Theorem \ref{thm0.1} (c)), this integral equals
$$\int\limits_{B_r\cap\{v>0\}}\nabla\bar v_{e_i}\cdot\nabla\(v_{e_i}*\psi_\epsilon\)
+\int\limits_{B_r\cap\{v<0\}}\nabla\bar
v_{e_i}\cdot\nabla\(v_{e_i}*\psi_\epsilon\).$$ On $S_v,\,\nabla
v=0,$ so that $\bar v_{e_i}$ will vanish on a neighborhood of $S_v$.
In fact, if $|v_{e_i}(x)|\geq\delta,\,z_0\in S_v$, then
$\delta\leq|v_{e_i}(x)-v_{e_i}(z_0)|\leq C|x-z_0|^\gamma$. In
$\F\setminus nbd(S_v),$ we have analyticity of $\F$ and a
well-defined normal, so that we can integrate by parts in the above
integrals, using Green's Theorem. We obtain for the above sum,
\begin{multline*}
-\int\limits_{B_r\cap\{v>0\}}\bar
v_{e_i}\Delta\(v_{e_i}*\psi_\epsilon\)
-\int\limits_{B_r\cap\{v<0\}}\bar
v_{e_i}\Delta\(v_{e_i}*\psi_\epsilon\) +\int\limits_{\p B_r}\bar
v_{e_i}\frac{\p}{\p\nu}\(v_{e_i}*\psi_\epsilon\)
\\+\int\limits_{B_r\cap\p\{v>0\}}\bar
v_{e_i}\frac{\p}{\p\nu}\(v_{e_i}*\psi_\epsilon\)
+\int\limits_{B_r\cap\p\{v<0\}}\bar
v_{e_i}\frac{\p}{\p\nu}\(v_{e_i}*\psi_\epsilon\).
\end{multline*}
The last two integrals cancel each other since the normals point in
opposite directions, in pieces of a real analytic surface. Thus, we
have obtained:
\begin{multline*}
\int_{B_r}\nabla\bar
v_{e_i}\cdot\nabla\(v_{e_i}*\psi_\epsilon\)=-\int\limits_{B_r\cap\{v>0\}}\bar
v_{e_i}\Delta\(v_{e_i}*\psi_\epsilon\)\\
-\int\limits_{B_r\cap\{v<0\}}\bar
v_{e_i}\Delta\(v_{e_i}*\psi_\epsilon\)+\int\limits_{\p B_r}\bar
v_{e_i}\frac{\p}{\p\nu}\(v_{e_i}*\psi_\epsilon\).
\end{multline*}
 We next average this identity in $r$, for
 $r\in\(\frac{r_0}2,\,\frac{3r_0}4\)$. We estimate first the
 averaged last term. Its absolute value is bounded by
 $$\ds c_nh\int\limits_{\frac{r_0}2\leq|x|\leq\frac{3r_0}4}|\nabla
 v_{e_i}*\psi_\epsilon|\leq c_nCh.$$
We next consider the absolute value of the averaged term of the
left-hand side, as $\epsilon\to0$. It converges to

$$\left|\frac4{r_0}\int_{r_0/2}^{3r_0/4}\int_{B_r}\nabla\bar v_{e_i}\cdot\nabla v_{e_i}\right|
\xrightarrow[\delta\to0]{}\left|\frac4{r_0}\int_{r_0/2}^{3r_0/4}\int_{B_r}\nabla\t
v_{e_i}\cdot\nabla v_{e_i}\right|$$ where
$$\t v_{e_i}=\begin{cases}
v_{e_i}& \text{if }|v_{e_i}|\leq h,\\
h & \text{otherwise.}
\end{cases}$$
This last expression is bounded below by
$c_n\int_{B_{r_0/2}}|\nabla\t v_{e_i}|^2$.

The absolute value of the sum on the averaged first two terms in the
right hand side converges (first letting $\epsilon\to0$ and then
$\delta\to0$) to

$$\left|\frac4{r_0}\int_{r_0/2}^{3r_0/4}\int_{B_r\cap\{v>0\}}\t v_{e_i}\Delta v_{e_i}
+\frac4{r_0}\int_{r_0/2}^{3r_0/4}\int_{B_r\cap\{v<0\}}\t
v_{e_i}\Delta v_{e_i}\right|.$$

But on $\{v>0\},\,\Delta v_{e_i}=\p_{e_i}f,$ on $\{v<0\},\,\Delta
v_{e_i}=-\p_{e_i}g.$ Hence, the above sum is bounded by
$$h\(\int_{B_{r_0}\cap\{v>0\}}|\Delta v_{e_i}|+\int_{B_{r_0}\cap\{v<0\}}|\Delta v_{e_i}|\)\leq Ch.$$
Finally, gathering terms and using that
$$\int_{B_{r_0/2}}|\nabla\t v_{e_i}|^2=\int_{B_{r_0/2}\cap\{|v_{e_i}|\leq h\}}|\nabla v_{e_i}|^2,$$
Claim \ref{claim2.45} follows.

We next complete the proof of the bound $H_{(n-1)}\(S_v^{(1)}\)\leq
C.$ Fix $z_0\in S_v^{(1)}$ and consider $S_v^{(1)}\cap
B\(z_0,\,r_0/4\)$. It suffices to prove our bound for this
intersection. For each $x_0$ in $S_v^{(1)}\cap B\(z_0,\,r_0/4\)$,
and each $0<r<r_0/100$, we consider the cover of $S_v^{(1)}\cap
B\(z_0,\,r_0/4\)$ by the balls $B(x_0,\,r).$ We can cover
$S_v^{(1)}\cap B\(z_0,\,r_0/4\)$ by finitely many such balls, and by
the Vitali covering Lemma, we can find $\t N$ disjoint balls
$B(x_i,\,r),\,x_i\in S_v^{(1)}\cap B\(z_0,\,r_0/4\)$ so that $\ds
S_v^{(1)}\cap B\(z_0,\,r_0/4\)\subset\bigcup_{i=1}^{\t
N}B(x_i,\,5r)$. The disjointness of $\{B(x_i,\,r)\}$ gives
$\ds\sum_{i=1}^{\t N}\chi_{B(x_i,5r)}(x)\leq c_n$. By Claim
\ref{claim2.44}, $|\nabla v(x)|\leq Cr$ in $B(x_i,\,5r)$. By
(\ref{eq2.2}) $|\Delta v|\geq C.$ We then have:
\begin{multline*}
c_n\t Ncr^n\leq\sum_i\int\limits_{B(x_i,5r)}\(\Delta v\)^2\leq
\int\limits_{\{|\nabla v(x)|\leq cr\}}\(\Delta v\)^2\sum_{i=1}^{\t
N}\chi_{B(x_i,5r)}\\\leq c_n\int\limits_{B(z_0,r_0/2)\cap\{|\nabla
v(x)|\leq cr\}}\(\Delta v\)^2\leq Cr,
\end{multline*}
by Claim \ref{claim2.45}. Thus, $\t Nr^{n-1}\leq C$, which gives our
Hausdorff measure bound.

To conclude this paper we give a simple result in the direction of
showing that better regularity results can hold for solutions of the
composite problem than for solutions of (\ref{eq2.2}) (see the end
of Remark \ref{rem2.40}). We will show that geometric assumptions on
$\Omega$ can ensure that for all solutions of the composite problem,
$S_u=\emptyset$ and thus $\F$ is real analytic and $u$ is $C^{1,1}$.

\begin{prop}\label{prop2.46}
Let $\Omega\subset\R^2$ have two axis of symmetry. Then for all
solutions $u$ of the composite problem (\ref{eq0.1}), (\ref{eq0.2})
we have $S_{u}=\emptyset$ and hence $\F$ is real analytic and $u\in
C^{1,1}.$
\end{prop}

\begin{proof}
We recall (see \cite{CGIKO}) that we say that $\Omega$ has an axis
of symmetry $L$ (which we take to be $\{x_1=0\}$) if whenever
$(x_1,\,x_2)$ belongs to $\Omega,$ so does $(-x_1,\,x_2)$ and the
set $\{x_1:(x_1,\,x_2)\in\Omega\}$ is either $\emptyset$ or an
interval $(-c,\,c)$ for each $x_2$. Let us give the proof, for
simplicity, in the case when the two axis $L_1,\,L_2$ are the $x_1$
and $x_2$ axis. It is shown in \cite{CGIKO}, Theorem 4, that any
solution $u$ is symmetric with respect to $x_1$ (and $x_2$) and $u$
is strictly decreasing in $x_1$, for $x_1\geq0$ (in $x_2$, for
$x_2\geq0$ ). (The strict decrease follows from $\alpha<\Lambda$,
see \cite{CGIKO}, the bottom of page 326). Because of the strict
decrease, $\frac{\p}{\p x_1}u(x_1,\,x_2)\neq0,\,x_1\neq0$ and
$\frac{\p}{\p x_2}u(x_1,\,x_2)\neq0$ for $x_2\neq0$. Thus, the only
possible point in $S_{u}$ is (0, 0). But, by the increase and
decrease described before $\ds u(0,\,0)=\sup_\Omega u$. Recall that
$D=\{0\leq u\leq c\},\,\F=\{u=c\}$. If $\ds c=\sup_{\Omega}
u,\,D=\Omega$, which contradicts $|D|=A<|\Omega|.$ Thus, (0,
0)$\notin\F$ and the Proposition follows.
\end{proof}

\setcounter{equation}{0}
\appendix

\section*{Appendix I}

The results (\ref{eqa1.9}), (\ref{eqa1.10}) are to be found in
\cite{CP}. They are reproduced here  for the reader's benefit.

We have the equation
\begin{equation}\label{eqa1.1}
-\Delta u_t +\alpha \chi_{D_t}u_t=\lambda(t)u_t \tag{A1.1}
\end{equation}
and the corresponding one for $u_0=u$, given by
\begin{equation}\label{eqa1.2}
-\Delta u +\alpha \chi_{D}u=\lambda u\tag{A1.2}
\end{equation}

where $\lambda(0)=\lambda$. We also note note that by our definition
of $D_t$,
\begin{equation}\label{eqa1.3}
\chi_{D_t}(x)=\chi_D(\phi_{-t}(x)).\tag{A1.3}
\end{equation}

We will set,
$$V(x)={{d\phi_t(x)}\over{dt}}\Big|_{t=0}$$
and assume  the vector field $V\in C^2(\Omega)$ and that $V$ is
supported in a compact set $S$.

Multiplying (\ref{eqa1.1}) by $u$, (\ref{eqa1.2}) by $u_t$ and
subtracting we get,
\begin{equation}\label{eqa1.4}
u_t\Delta u-u\Delta
u_t+\alpha(\chi_D(\phi_{-t}(x))-\chi_D(x))uu_t=(\lambda(t)-\lambda)uu_t.\tag{A1.4}
\end{equation}
 We integrate \ref{eqa1.4} over $\Omega$.  Since, $u=u_t=0$ on
$\partial\Omega$ we get
$$\int_\Omega u_t\Delta u-u\Delta u_t=0$$
Thus the integral over $\Omega$ of (\ref{eqa1.4}) becomes,
\begin{equation}\label{eqa1.5}
\int_\Omega \alpha(\chi_D(\phi_{-t}(x))-\chi_D(x))uu_t=
(\lambda(t)-\lambda)\int_\Omega uu_t.\tag{A1.5}
\end{equation}
Now from (\ref{eqa1.1}) we notice that if we normalize our functions
$\|u_t\|_2=1$ as we certainly can, we always have $\|u_t\|_{2,2}\leq
C$. Now,
$$\left|\int_\Omega (uu_t-u^2)\right|\leq \int_\Omega u|u-u_t|$$
In a
tubular neighborhood of $\partial\Omega$, ${\cal U}$ we have,
$$\int_{\cal U} u|u-u_t|\leq C\(\int_{\cal U} u^2\)^{1/2}\leq \epsilon$$
Outside ${\cal U}$ by the uniform $W^{2,2}$ bounds of $u_t$ we have
strong convergence of $u_t$ to $u$ in $L^2$. Thus we have,
\begin{equation}\label{eqa1.6}
\lim_{t\to 0}\int_\Omega uu_t=\int_\Omega u^2=1.\tag{A1.6}
\end{equation}

Now we change variables in the left side of (\ref{eqa1.5}). We set
$\phi_{-t}(x)=y$. Thus, $x=\phi_{-t}^{-1}(y)$. Thus the left side of
(\ref{eqa1.5}) becomes,
$$\alpha\int_\Omega \chi_D(x)( h_t(\phi_{-t}^{-1}(x))J_t(x)-h_t(x))dx.$$
Here we have set $h_t(x)=uu_t$ and $J_t(x)$ is the determinant of
the Jacobian matrix of the transformation $y=\phi_{-t}^{-1}(x)$.
Since $\phi_0(x)=I$ the identity, it is well-known that,
\begin{equation}\label{eqa1.7}
J_t(x)=1+t\, {\rm div}\, V+O\(t^2\).\tag{A1.7}\end{equation} See for
example Lemma 1(pg 69) in \cite{A}, in fact (\ref{eqa1.7}) is an
elementary consequence of the fact that for a $n\times n$ matrix
$B$, $\det (I-tB)^{-1}=1+t\, {\rm trace}\ B+O\(t^2\)$. Since $h_t\in
C^{1,\beta}$, we see that,
$$h_t\(\phi_{-t}^{-1}(x)\)J_t(x)-h_t(x)= t\((V\cdot \nabla)h_t(x)+h_t(x){\rm div}\, V \)+o(t)$$
Thus on division by $t$ and letting $t\to 0$ we see easily,
$$ \lim_{t\to0}{{h_t\(\phi_{-t}^{-1}(x)\)J_t(x)-h_t(x)}\over{t}}=(V\cdot\nabla)
(u^2)+u^2{\rm div}\, V.$$
The term on the right above is,
$${\rm div}\ \(Vu^2\).$$
Thus dividing (\ref{eqa1.5}) by $t$ and using (\ref{eqa1.6}) we
easily get,
$$\lambda^\prime(0)=\alpha \int_D{\rm div}\, \(Vu^2\)=\alpha\int_{D\cap
S}{\rm div}\, \(Vu^2\). $$ By the hypothesis that the part of the
boundary of $\partial D$ that lies inside the support of $V$ is
regular enough to have a bonafide unit outer normal $\nu$, and
Green's theorem, the last integral above yields,
\begin{equation}\label{eqa1.8}
\lambda'(0)=\alpha\int_{S\cap
\partial D}\langle V,\nu\rangle u^2.\tag{A1.8}\end{equation}

Now consider,
$$|D_t|-|D|=\int_\Omega (\chi_D(\phi_{{-t}}(x))-\chi_D(x))\ dx$$
Change variables in the integral above as before to get,
$$\int_D\(J_t(x)-1\)\ dx.$$
By (\ref{eqa1.7}) again we see the integral above is,
   $$t\int_D{\rm div}\, V \ dx+O\(t^2\).$$
Thus we easily get,
\begin{equation}\label{eqa1.9}
{{d}\over{dt}}(|D_t|)|_{t=0}=\int_D {\rm div}\, V\ dx=\int_{S\cap
\partial D}\langle V,\nu\rangle\, d\sigma.\tag{A1.9}\end{equation}

If $u=c$ along $\partial D$, combining (\ref{eqa1.8}) and
(\ref{eqa1.9}) we get,
\begin{equation}\label{eqa1.10}
\lambda'(0)=\alpha c^2\frac d{dt}\(|D_t|\)|_{t=0}=\alpha c^2\int_{\p
D}\langle V,\,\nu\rangle\ d\sigma.\tag{A1.10}
\end{equation}

\section*{Appendix II}

We use Blank's \cite{B} notation. We have,
$$f_1(\theta)= C_+\sin (2\theta+D_+)+\gamma,\ f_1>0$$
and also,
$$f_2(\theta)=C_-\sin (2\theta+D_-)+\mu,\ f_2<0$$
Now we focus on the interval $[0, 2\pi/3]$. First for $\theta_0\in
(0,2\pi/3)$, we know $f_1(\theta_0)=f_2(\theta_0)=0$ and
$f_1^\prime(\theta_0)=f_2^\prime(\theta_0)=0$. We get,
$$ C_+\sin (2\theta_0+D_+)+\gamma=C_-\sin (2\theta_0+D_-)+\mu=0$$
and,
$$C_+\cos(2\theta_0+D_+)=C_-\cos(2\theta_0+D_-)$$
This leads after squaring and adding both equations to,
\begin{equation}\label{eqa2.1}
C_+^2-C_-^2=\gamma^2-\mu^2.\tag{A2.1}\end{equation} Next because
$f_1(0)=f_1(\theta_0)=0$, we get,
\begin{equation}\label{eqa2.2}
C_+\sin(D_+)=-\gamma,\ D_+=\arcsin
(-\gamma/C_+)\tag{A2.2}\end{equation} and also we have,
\begin{equation}\label{eqa2.3}\theta_0=\pi/2+\arcsin(\gamma/C_+).\tag{A2.3}\end{equation}
Since $f_2(\theta_0)=0$, inserting the value of $\theta_0$ from
(\ref{eqa2.3}3) in the expression for $f_2$, we see,
\begin{equation}\label{eqa2.4}
D_-=\arcsin(\mu/C_-)-2\arcsin(\gamma/C_+).\tag{A2.4}\end{equation}
Now lastly $f_2(2\pi/3)=0$, so,
$$C_-\sin(4\pi/3+D_-)=-\mu$$
We get,
\begin{equation}\label{eqa2.5}C_-\ =\
\mu/\sin(\pi/3+D_-)\tag{A2.5}\end{equation} Now assume
$|C_+|>10^6\(|\gamma|+|\mu|\)$, then by (\ref{eqa2.1}),
$|C_-|>10^6(|\gamma|+|\mu|)$. Thus, from (\ref{eqa2.4}), $|D_-|\leq
\pi/20$. From (\ref{eqa2.5}) we get,
$$|C_-|\leq 2|\mu|$$
And we get a bound on $|C_+|$ from (\ref{eqa2.1}) again.

\bibliographystyle{amsalpha}
\bibliography{bibweak}

\newcommand{\etalchar}[1]{$^{#1}$}
\providecommand{\bysame}{\leavevmode\hbox to3em{\hrulefill}\thinspace}
\providecommand{\MR}{\relax\ifhmode\unskip\space\fi MR }
\providecommand{\MRhref}[2]{%
  \href{http://www.ams.org/mathscinet-getitem?mr=#1}{#2}
}
\providecommand{\href}[2]{#2}
\begin{thebibliography}{CGI{\etalchar{+}}00}

\bibitem[Arn97]{A}
V.~I. Arnold, \emph{Mathematical methods of classical mechanics (graduate texts
  in mathematics)}, Springer, September 1997.

\bibitem[AW05]{AW}
J.~Andersson and G.~S. Weiss, \emph{Cross-shaped and degenerate singularities
  in an unstable elliptic free boundary problem}, Preprint,
  http://arxiv.org/abs/math.AP/0508184, 2005.

\bibitem[Bla04]{B}
I.~Blank, \emph{Eliminating mixed asymptotics in obstacle type free boundary
  problems}, Comm. PDE \textbf{29} (2004), no.~7-8, 1167--1186.

\bibitem[CGI{\etalchar{+}}00]{CGIKO}
S.~Chanillo, D.~Grieser, M.~Imai, K.~Kurata, and I.~Ohnishi, \emph{Symmetry
  breaking and other phenomena in the optimization of eigenvalues for composite
  membranes}, Comm. Math. Phys. \textbf{214} (2000), 315--337.

\bibitem[CGK00]{CGK}
S.~Chanillo, D.~Grieser, and K.~Kurata, \emph{The free boundary problem in the
  optimization of composite membranes}, Contemp. Math. \textbf{268} (2000),
  61--81.

\bibitem[CP]{CP}
S.~Chanillo and R.~Pedroza, \emph{Hadamard's formulae for composite membranes},
  Preprint.

\bibitem[Giu84]{G}
E.~Giusti., \emph{Minimal surfaces and functions of bounded variation},
  Monographs in Mathematics, Birkhauser, 1984.

\bibitem[GT83]{GT}
D.~Gilbarg and N.~Trudinger, \emph{Elliptic partial differential equations of
  second order}, Springer-Verlag, 1983.

\bibitem[Kat73]{K}
T~Kato, \emph{Schrödinger operators with singular potentials}, Proceedings of
  the International Symposium on Partial Differential Equations and the
  Geometry of Normed Linear Spaces (Jerusalem, 1972), vol.~13, 1973,
  pp.~135--148.

\bibitem[MW05]{MW}
R.~Monneau and G.~S. Weiss, \emph{An unstable elliptic free boundary problem
  arising in solid combustion}, Preprint, http://arxiv.org/abs/math.AP/0507315,
  2005.

\bibitem[Sha]{S}
H.~Shahgholian, \emph{The singular set for the composite membrane problem},
  Preprint.

\bibitem[Wei98]{W}
G.~S. Weiss, \emph{Partial regularity for weak solutions of an elliptic free
  boundary problem}, Comm. PDE \textbf{23} (1998), no.~3-4, 439--455.

\end{thebibliography}

\end{document}